\documentclass[12pt]{article}

\usepackage{amsfonts}
\usepackage{amsmath}
\usepackage{amsthm}

\newcommand{\supp}{\mathrm{supp}}
\numberwithin{equation}{section}
\renewcommand{\O}{\mathcal{O}}
\newtheorem{theorem}{Theorem}[section]
\newtheorem{proposition}{Proposition}[section]
\DeclareMathOperator{\diag}{diag}
\def\ds{\displaystyle}

\begin{document}

\title{Ladder operators and differential equations
for multiple orthogonal polynomials}
\author{Galina Filipuk\footnotemark[1],\quad Walter Van Assche\footnotemark[2],
\quad Lun Zhang\footnotemark[2]}
\date{\today}

\maketitle
\renewcommand{\thefootnote}{\fnsymbol{footnote}}
\footnotetext[1]{Faculty of Mathematics, Informatics and Mechanics,
University of Warsaw, Banacha 2, Warsaw, 02-097, Poland. E-mail:
filipuk@mimuw.edu.pl} \footnotetext[2]{Department of Mathematics,
KU Leuven, Celestijnenlaan 200B, B-3001 Leuven,
Belgium. E-mail: \{Walter.VanAssche, lun.zhang\}@wis.kuleuven.be}

\begin{abstract}
In this paper, we obtain the ladder operators and associated
compatibility conditions for the type I and the type II multiple
orthogonal polynomials. These ladder equations extend known results
for orthogonal polynomials and can be used to derive the
differential equations satisfied by multiple orthogonal polynomials.
Our approach is based on Riemann-Hilbert problems and the
Christoffel-Darboux formula for multiple orthogonal polynomials, and
the nearest-neighbor recurrence relations. As an illustration, we
give several explicit examples involving multiple Hermite and
Laguerre polynomials, and multiple orthogonal polynomials with
exponential weights and cubic potentials.

\vspace{2mm} \textbf{Keywords:} Ladder operators, differential
equations, multiple orthogonal polynomials, Riemann-Hilbert
problems, nearest-neighbor recurrence relations.

\end{abstract}

\section{Introduction }

%----------------------------------------------------------------
\subsection{Multiple orthogonal polynomials}\label{sec:MOP intro}

Multiple orthogonal polynomials are polynomials of one variable
which are defined by orthogonality relations with respect to $r$
different measures $\mu_1,\mu_2,\allowbreak \ldots,\mu_r$, where
$r \geq 1$ is a positive integer. As a generalization of
orthogonal polynomials \cite{Chi,Gau,Sze}, multiple orthogonal
polynomials originated from Hermite-Pad\'e approximation in the context
of irrationality and transcendence proofs
in number theory. They were further developed in approximation
theory, cf.
\cite{Apt98,Apt99,Bust,Gon,NikSor} and surveys
\cite{Apt,WVA99,WVA06}. During the past few years, multiple
orthogonal polynomials have also arisen in a natural way in
certain models from mathematical physics, including random matrix
theory, non-intersecting paths, etc.; we refer to
\cite{Kuijl,Kuijl2} and references therein for the progress of
this subject.%aspect.

Let $\vec{n} = (n_1,n_2,\ldots,n_r) \in \mathbb{N}^r$ be a
multi-index of size $|\vec{n}| = n_1+n_2+\ldots+n_r$ and suppose
$\mu_1,\mu_2,\ldots,\mu_r$ are $r$ measures with supports on certain
simple curves in the complex plane. There are two types of multiple
orthogonal polynomials. The \emph{type I multiple orthogonal
polynomials} are given by the vector $(A_{\vec{n},1},\ldots,
A_{\vec{n},r})$, where $A_{\vec{n},j}$ is a polynomial of degree
$\leq n_j-1$, for which
\begin{equation}  \label{eq:1.6}
     \sum_{j=1}^r \int  x^k  A_{\vec{n},j}(x)   \, d\mu_j(x) = 0,
    \qquad k=0,1,\ldots,|\vec{n}|-2.
\end{equation}
We use the normalization
\begin{equation} \label{eq:1.7}
       \sum_{j=1}^r \int  x^{|\vec{n}|-1}  A_{\vec{n},j}(x)  \, d\mu_j(x) = 1.
\end{equation}
The equations \eqref{eq:1.6}--\eqref{eq:1.7} are a linear system of $|\vec{n}|$ equations for the unknown
coefficients of $A_{\vec{n},1},\ldots,A_{\vec{n},r}$. This system has a unique solution if the matrix of mixed moments
\begin{equation*} \label{eq:1.8}
    M_{\vec{n}} = \left( \begin{matrix}
                   \nu_0^{(1)} & \nu_1^{(1)} & \cdots & \nu_{n_1-1}^{(1)} \\
                   \nu_1^{(1)} & \nu_2^{(1)} & \cdots & \nu_{n_1}^{(1)} \\
                \vdots & \vdots & \cdots & \vdots \\
                   \nu_{|\vec{n}|-1}^{(1)} & \nu_{|\vec{n}|}^{(1)} & \cdots & \nu_{|\vec{n}|+n_1-2}^{(1)}
                  \end{matrix} \right|
                % \left. \begin{matrix}
                %   \nu_0^{(2)} & \nu_1^{(2)} & \cdots & \nu_{n_2-1}^{(2)} \\
                %   \nu_1^{(2)} & \nu_2^{(2)} & \cdots & \nu_{n_2}^{(2)} \\
            %    \vdots & \vdots & \cdots & \vdots \\
                %   \nu_{|\vec{n}|-1}^{(2)} & \nu_{|\vec{n}|}^{(1)} & \cdots & \nu_{\vec{n}+n_2-2}^{(2)}
                %  \end{matrix} \right|
        \cdots \cdots
                 \left| \begin{matrix}
                   \nu_0^{(r)} & \nu_1^{(r)} & \cdots & \nu_{n_r-1}^{(r)} \\
                   \nu_1^{(r)} & \nu_2^{(r)} & \cdots & \nu_{n_r}^{(r)} \\
                \vdots & \vdots & \cdots & \vdots \\
                   \nu_{|\vec{n}|-1}^{(r)} & \nu_{|\vec{n}|}^{(r)} & \cdots & \nu_{|\vec{n}|+n_r-2}^{(r)}
                  \end{matrix} \right) ,
\end{equation*}
where
\[    \nu_n^{(j)} = \int x^n\, d\mu_j(x), \qquad j=1,\ldots,r,  \]
is not singular, in which case we call the multi-index $\vec{n}$ a
\textit{normal index}. The \emph{type II multiple orthogonal
polynomial} is the monic polynomial $P_{\vec{n}}(x) = x^{|\vec{n}|}
+ \cdots$ of degree $|\vec{n}|$ for which
\begin{equation} \label{eq:orthogonality of Pn}
\begin{aligned}
    \int P_{\vec{n}}(x) x^k\, d\mu_1(x) &= 0, \qquad k=0,1,\ldots,n_1-1,  \\
    \vdots  & & \\
    \int P_{\vec{n}}(x) x^k\, d\mu_r(x) &= 0, \qquad k=0,1,\ldots,n_r-1.
    \end{aligned}
\end{equation}
This gives a linear system of $|\vec{n}|$ equations for the
$|\vec{n}|$ unknown coefficients of $P_{\vec{n}}$ and the matrix of
this linear system is the transpose $M_{\vec{n}}^T$ of
\eqref{eq:1.8}, hence $P_{\vec{n}}$ exists and is unique whenever
$\vec{n}$ is a normal index.

Suppose that the $r$ measures $\mu_1,\ldots,\mu_r$ are all
absolutely continuous with respect to a measure $\mu$ and that
$d\mu_j(x) = w_j(x)\, d\mu(x)$. The type I and II multiple
orthogonal polynomials satisfy the following biorthogonality:
\begin{equation}\label{eq:biorthogonality}
\int P_{\vec{n}}(x) Q_{\vec{m}}(x)\, d\mu(x)
  = \begin{cases} 0, & \quad\textrm{if $\vec{m} \leq \vec{n}$,} \\
                  0, & \quad \textrm{if $|\vec{n}| \leq |\vec{m}|-2$,} \\
                  1, & \quad \textrm{if $|\vec{m}|=|\vec{n}|+1$,}
    \end{cases}
\end{equation}
where
\begin{equation}\label{def:Qn}
Q_{\vec{n}}(x) = \sum_{j=1}^r A_{\vec{n},j}(x)w_j(x),
\end{equation}
see \cite[Theorem 23.1.6]{Ismail}.

For more information about multiple orthogonal polynomials, we refer
to Aptekarev et al. \cite{Apt,AptBraWVA}, Coussement and Van Assche
\cite{WVAEC}, Nikishin and Sorokin \cite[Chapter 4, \S 3]{NikSor},
and Ismail \cite[Chapter 23]{Ismail}. Throughout this paper, we
shall assume that all multi-indices are normal.

\subsection{Recurrence relations }\label{sec:recurrence relation}

There are several recurrence relations for the type I and the type
II multiple orthogonal polynomials. For type II, we have the
following nearest-neighbor recurrence relations
\begin{equation}\label{eq:1.10}
\begin{aligned}
   x P_{\vec{n}}(x) &= P_{\vec{n}+\vec{e}_1}(x) + b_{\vec{n},1}P_{\vec{n}}(x) + \sum_{j=1}^r a_{\vec{n},j} P_{\vec{n}-\vec{e}_j}(x),
   \\
      \vdots & &  \\
   x P_{\vec{n}}(x) & =  P_{\vec{n}+\vec{e}_r}(x) + b_{\vec{n},r}P_{\vec{n}}(x) + \sum_{j=1}^r a_{\vec{n},j} P_{\vec{n}-\vec{e}_j}(x),
\end{aligned}
\end{equation}
where $\vec{e}_j = (0,\ldots,0,1,0,\ldots,0)$ is the $j$-th standard
unit vector with 1 on the $j$-th entry,
$(a_{\vec{n},1},\ldots,a_{\vec{n},r})$ and
$(b_{\vec{n},1},\ldots,b_{\vec{n},r})$ are the recurrence
coefficients. These recurrence relations were derived in
\cite[Theorem 23.1.11]{Ismail} where it was shown that
\begin{equation} \label{eq:a(n,j)}
     a_{\vec{n},j} = \frac{\displaystyle \int x^{n_j} P_{\vec{n}}(x)\, d\mu_j(x)}
                            {\displaystyle \int x^{n_j-1} P_{\vec{n}-\vec{e}_j}(x)\, d\mu_j(x)},
\end{equation}
and
\begin{equation} \label{eq:1.12}
    b_{\vec{n},j} = \int xP_{\vec{n}}(x) Q_{\vec{n}+\vec{e}_j}(x)\, d\mu(x),
\end{equation}
for $j=1,\ldots,r$, where $Q_{\vec{n}}$ is given in \eqref{def:Qn}.

We have similar recurrence relations for the type I multiple
orthogonal polynomials:
\begin{equation} \label{eq:1.10Q}
\begin{aligned}
   x Q_{\vec{n}}(x) &= Q_{\vec{n}-\vec{e}_1}(x) + b_{\vec{n}-\vec{e}_1,1}Q_{\vec{n}}(x) + \sum_{j=1}^r a_{\vec{n},j} Q_{\vec{n}+\vec{e}_j}(x),
    \\
      \vdots & &  \\
   x Q_{\vec{n}}(x) & =
    Q_{\vec{n}-\vec{e}_r}(x) + b_{\vec{n}-\vec{e}_r,r}Q_{\vec{n}}(x) +
    \sum_{j=1}^r a_{\vec{n},j} Q_{\vec{n}+\vec{e}_j}(x).
\end{aligned}
\end{equation}
Observe that the same recurrence coefficients
$(a_{\vec{n},1},\ldots,a_{\vec{n},r})$ are used but that there is a
shift in the other recurrence coefficients $(b_{\vec{n}-\vec{e}_1,1},
\ldots, b_{\vec{n}-\vec{e}_r,r})$.

If we multiply (\ref{eq:1.10}) by $Q_{\vec{n}}(x)$ and integrate,
then the biorthogonality \eqref{eq:biorthogonality} gives
\begin{equation*}  \label{eq:suma}
   \int x P_{\vec{n}}(x) Q_{\vec{n}}(x)\, d\mu(x) = \sum_{j=1}^r a_{\vec{n},j} .
\end{equation*}
The orthogonality properties of $P_{\vec{n}}$ in
\eqref{eq:orthogonality of Pn} imply that
\[    \int P_{\vec{n}}(x)Q_{\vec{n}+\vec{e}_j}(x)\, d\mu(x)
= \int P_{\vec{n}}(x) A_{\vec{n}+\vec{e}_j,j}(x)\, d\mu_j(x), \]
hence, by \eqref{eq:biorthogonality},
\begin{equation}  \label{eq:kappa}
    1 = \kappa_{\vec{n}+\vec{e}_j,j} \int  x^{n_j} P_{\vec{n}}(x)\,
    d\mu_j(x),
\end{equation}
where $A_{\vec{n},j}(x) = \kappa_{\vec{n},j} x^{n_j-1} + \cdots$.
This, together with \eqref{eq:a(n,j)}, provides an alternative
representation of $a_{\vec{n},j}$:
\begin{equation}  \label{eq:an representation}
   a_{\vec{n},j} = \frac{ \int x^{n_j} P_{\vec{n}}(x)\, d\mu_j(x) }
   { \int x^{n_j-1} P_{\vec{n}-\vec{e}_j}(x)\, d\mu_j(x) }
   = \frac{\kappa_{\vec{n},j}}{\kappa_{\vec{n}+\vec{e}_j,j}}.
\end{equation}

\subsection{Ladder equations for orthogonal polynomials}
Given a single positive measure $w(x)dx$, it is well-known that
there exist monic orthogonal polynomials $P_n(x)$ of degree $n$ in
$x$ such that
\begin{equation} \label{orthogonality monic}
\int P_m(x) P_n(x) w(x)\, dx = h_n \delta_{m,n}, \quad h_n>0, \quad
m,n=0,1,2,\ldots.
\end{equation}
They are characterized by the three-term recurrence relation:
\begin{equation} \label{recurrence}
x P_n(x) = P_{n+1}(x) + \alpha_n P_n(x) + \beta_n P_{n-1}(x),
\end{equation}
where
\begin{equation}\label{recu:monic}
\alpha_n=\frac{1}{h_n}\int x P_n^2(x) w(x)\, dx,\qquad
\beta_n=\frac{1}{h_{n-1}}\int x P_n(x)P_{n-1}(x) w(x)\, dx,
\end{equation}
and the initial condition is taken to be $\beta_0 P_{-1}(x) :=0$.
Suppose that $w$ vanishes at the end points of the orthogonality
interval, then it was shown in \cite{chen1} that $P_n$ satisfy the
following ladder equations:
\begin{align}
\left( \frac{d}{dx} + B_n(x) \right) P_n(x) & = \beta_n A_n(x)
P_{n-1}(x), \label{ladder1} \\ \left( \frac{d}{dx} - B_n(x) -
\textsf{v}'(x) \right) P_{n-1}(x) & = - A_{n-1}(x) P_n(x)
\label{ladder2}
\end{align}
with $\textsf{v}(x):=-\ln w(x)$  and  \begin{align} A_n(x) & :=
\frac{1}{h_n} \int \frac{\textsf{v}'(x) - \textsf{v}'(y)}{x-y} \
[P_n(y)]^2 w(y) dy, \label{an-def}\\ B_n(x) & := \frac{1}{h_{n-1}}
\int \frac{\textsf{v}'(x) - \textsf{v}'(y)}{x-y} \ P_{n-1}(y) P_n(y)
w(y) dy. \label{bn-def}
\end{align}
Furthermore, the functions $A_n$ and $B_n$ defined  by
(\ref{an-def}) and
(\ref{bn-def}) satisfy % the following conditions
$$
B_{n+1}(z) + B_n(z)  = (z- \alpha_n) A_n(z) -
\textup{\textsf{v}}\,'(z), \eqno(S_1)
$$
$$
1+ (z- \alpha_n) [B_{n+1}(z) - B_n(z)] = \beta_{n+1} A_{n+1}(z) -
\beta_n A_{n-1}(z).\eqno(S_2)
$$
The conditions $S_1$ and $S_2$ are usually called the compatibility
conditions for the ladder equations. Similar relations were also
achieved for discrete orthogonal polynomials \cite{INS},
$q$-orthogonal polynomials \cite{chen+ismail08} and matrix
orthogonal polynomials \cite{DI,Gru}.

The motivation of deriving such ladder equations is two-fold. On the
one hand, the differential recurrence relations \eqref{ladder1} and
\eqref{ladder2}, together with recurrence relation
\eqref{recurrence}, will lead to a linear differential equation
satisfied by $P_n$. On the other hand, ladder operators has been
successfully applied to establish the connections between
Painlev\'e equations and recurrence coefficients of certain
orthogonal polynomials; cf.
\cite{BCK,BFV,ChenIts,chen+zhang,Dai+Zhang,GF+Walter1,GF+Walter2,GVZ}
for recent applications.

It is the aim of this paper to find and explore ladder relations for
multiple orthogonal polynomials. We hope our results will be helpful
in the further study of recurrence coefficients of multiple orthogonal
polynomials and its relation to certain integrable systems. In the
literature, the first attempt in this direction was done by
Coussement and Van Assche \cite{JCWVA,WVAEC}, where they derived
raising and lowering operators for type II multiple orthogonal
polynomials with respect to classical weights. The results there,
however, are different from ours, since they consider relations
between multiple orthogonal polynomials associated with different
weights in the same class.

The results of this paper will be presented in the next section.

\section{Statement of the results}

\subsection{Ladder equations for multiple orthogonal polynomials}
\label{sec:ladder equations}

Our first result concerns the ladder equations for the type II
multiple orthogonal polynomials.

\begin{theorem}\label{thm:ladder operators for type II}
Let $\mu_1,\ldots,\mu_r$ be $r$ measures that are absolutely
continuous with weights $w_1,\ldots,w_r$ and each $w_i$ vanishes at
the endpoints of the support of $\mu_i$. Suppose that all the
indices $\vec{n}=(n_1,\ldots,n_r)\in\mathbb{N}^r$ are normal and the
functions
\[  \{ w_1,xw_1,\ldots,x^{n_1-1}w_1,w_2,xw_2,\ldots,x^{n_2-1}w_2,\ldots,
     w_r,xw_r,\ldots,x^{n_r-1}w_r \} \]
are linearly independent, then  we have the following lowering
equation for the type II multiple orthogonal polynomials:
\begin{align}\label{lowering II}
P_{\vec{n}}'(x)=&~P_{\vec{n}}(x)\int
P_{\vec{n}}(t) \sum_{k=1}^r
A_{\vec{n},k}(t)\frac{v_k'(t)-v_k'(x)}{x-t} w_k(t)\,dt
\nonumber
\\ &-\sum_{j=1}^r a_{\vec{n},j} P_{\vec{n}-\vec{e}_j}(x)\int
P_{\vec{n}}(t)   \sum_{k=1}^r
A_{\vec{n}+\vec{e}_j,k}(t)\frac{v_k'(t)-v_k'(x)}{x-t}
 w_k(t)\,dt,
\end{align}
where $v_k(x):=-\ln w_k(x)$ and $a_{\vec{n},j}$ are the recurrence
coefficients given in \eqref{eq:a(n,j)}.

The following raising equations for the type II
multiple orthogonal polynomials hold: for $i=1\ldots,r$ one has
\begin{align}\label{raising II}
P_{\vec{n}-\vec{e}_i}'(x)=&~P_{\vec{n}}(x)\int
P_{\vec{n}-\vec{e}_i}(t) \sum_{k=1}^r
A_{\vec{n},k}(t)\frac{v_k'(t)-v_k'(x)}{x-t} w_k(t)\,dt\nonumber\\
&-\sum_{j=1}^r \Big(a_{\vec{n},j}\int P_{\vec{n}-\vec{e}_i}(t)
  \sum_{k=1}^r
A_{\vec{n}+\vec{e}_j,k}(t)\frac{v_k'(t)-v_k'(x)}{x-t} w_k(t)\,dt\nonumber
\\ &~~~~~~~~~~~~~~~~~~-v_i'(x)\delta_{i,j} \Big) P_{\vec{n}-\vec{e}_j}(x),
\end{align}
where $\delta_{i,j}$ is the Kronecker delta.
\end{theorem}

Although our assumption on the weight functions may seem a little bit
involved, they include the so-called Angelesco systems
and AT systems that are normal for all indices; cf. \cite[Sections
23.1.1--23.1.2]{Ismail} for an introduction of these two weight
systems. Moreover, in the case of $r=1$, suppose
$d\mu_1(x)=w(x)dx$, we have, in the notations of multiple
orthogonal polynomials introduced in Sections \ref{sec:MOP intro}
and \ref{sec:recurrence relation},
\begin{equation*}
P_{\vec{n}}(x)=P_n(x),\quad
A_{\vec{n},1}(x)=\frac{P_{n-1}(x)}{h_{n-1}}, \quad
b_{\vec{n},1}=\alpha_n, \quad a_{\vec{n},1}=\beta_n,
\end{equation*}
where $P_n$, $h_n$, $\alpha_n$ and $\beta_n$ are given in
\eqref{orthogonality monic} and \eqref{recu:monic}. Inserting the
above representations into \eqref{lowering II} and \eqref{raising
II}, we recover equations \eqref{ladder1} and \eqref{ladder2} by
noting the fact that $\beta_n=h_n/h_{n-1}$.

The results in Theorem \ref{thm:ladder operators for type II} can be
summarized in vector form. To this end, we introduce an
$(r+1)\times(r+1)$ matrix
\begin{equation}\label{def:N}
N(x):=N(\vec{n};x)=(N_{ij}(x))_{0\leq i,j\leq r},
\end{equation}
where
\begin{equation}\label{def:Nij}
N_{ij}(x)=\left\{
            \begin{array}{ll}
              \ds \int P_{\vec{n}-\vec{e}_i}(t)\sum\limits_{k=1}^r A_{\vec{n},k}(t)
               \frac{v_k'(t)-v_k'(x)}{x-t}w_k(t)\,dt, & \hbox{$j=0$,}
\\
             \ds a_{\vec{n},j}\int P_{\vec{n}-\vec{e}_i}(t) \sum\limits_{k=1}^r
              A_{\vec{n}+\vec{e}_j,k}(t)\frac{v_k'(t)-v_k'(x)}{t-x}w_k(t)\,dt
              +v_i'(x)\delta_{i,j}, & \hbox{$j\neq 0$.}
            \end{array}
          \right.
\end{equation}
Here, it is understood that $\vec{e}_0=\vec{0}$. We then have from
\eqref{lowering II} and \eqref{raising II} that
\begin{equation}\label{eq:matrix form type II}
\textbf{P}_{\vec{n}}'(x) =N(x) \textbf{P}_{\vec{n}}(x),
\end{equation}
where
\begin{equation}\label{def:vec P}
\textbf{P}_{\vec{n}}(x):=\left(P_{\vec{n}}(x),P_{\vec{n}-\vec{e}_1}(x),\ldots,
P_{\vec{n}-\vec{e}_r}(x)\right)^T.
\end{equation}

We have similar results for the type I multiple orthogonal
polynomials.

\begin{theorem}\label{thm:ladder operators for type I}
Let $\mu_1,\ldots,\mu_r$ be $r$ measures as given in Theorem
\ref{thm:ladder operators for type II}. For $l=1,\ldots,r$, the
associated type I multiple orthogonal polynomials satisfy the
following raising equations:
\begin{align}\label{raising I}
A_{\vec{n},l}'(x)=&~-A_{\vec{n},l}(x)\int
P_{\vec{n}}(t) \sum_{k=1}^r
A_{\vec{n},k}(t)\frac{v_k'(t)-v_k'(x)}{x-t}w_k(t)\,dt
\nonumber
\\ &+\sum_{j=1}^r a_{\vec{n},j} A_{\vec{n}+\vec{e}_j,l}(x)\int
P_{\vec{n}-\vec{e}_j}(t)  \sum_{k=1}^r
A_{\vec{n},k}(t)\frac{v_k'(t)-v_k'(x)}{x-t} w_k(t)\,dt,
\end{align}
where $v_k(x):=-\ln w_k(x)$ and $a_{\vec{n},j}$ are the recurrence
coefficients given in \eqref{eq:a(n,j)}.

The lowering equations for the type I multiple orthogonal
polynomials are given by
\begin{align}\label{lowering I}
A_{\vec{n}+\vec{e}_i,l}'(x)=&~-A_{\vec{n},l}(x)\int
P_{\vec{n}}(t) \sum_{k=1}^r
A_{\vec{n}+\vec{e}_i,k}(t)\frac{v_k'(t)-v_k'(x)}{x-t} w_k(t)\,dt\nonumber\\
&+\sum_{j=1}^r\Big(a_{\vec{n},j}\int P_{\vec{n}-\vec{e}_j}(t)
  \sum_{k=1}^r
A_{\vec{n}+\vec{e}_i,k}(t)\frac{v_k'(t)-v_k'(x)}{x-t} w_k(t)\,dt\nonumber
\\ &~~~~~~~~~~~~~~~~~~~~-v_i'(x)\delta_{i,j} \Big)A_{\vec{n}+\vec{e}_j,l}(x),
\end{align}
for $i=1,\ldots,r$.
\end{theorem}

In vector form, Theorem \ref{thm:ladder operators for type I}
reads
\begin{equation}\label{eq:matrix form type I}
\frac{d}{dx}
 \begin{pmatrix} -A_{\vec{n},l}(x) \\
a_{\vec{n},1} A_{\vec{n}+\vec{e}_1,l}(x)\\
\vdots \\
a_{\vec{n},r} A_{\vec{n}+\vec{e}_r,l}(x)
\end{pmatrix}
= N^T(x)\begin{pmatrix} A_{\vec{n},l}(x) \\
-a_{\vec{n},1} A_{\vec{n}+\vec{e}_1,l}(x)\\
\vdots \\
-a_{\vec{n},r} A_{\vec{n}+\vec{e}_r,l}(x)
\end{pmatrix}, \qquad l=1,\ldots,r.
\end{equation}
%where
%\begin{equation}\label{def:tilde N}
%\widetilde N(x)=\diag\left(-1, \frac{1}{a_{\vec{n},1}}, \ldots,
%\frac{1}{ a_{\vec{n},r}}\right )N^{T}(x) \diag(1, -a_{\vec{n},1},
%\ldots, - a_{\vec{n},r})
%\end{equation}
%with $N$ being defined in \eqref{def:N}--\eqref{def:Nij}.

 From the lowering and raising equations stated in Theorems
\ref{thm:ladder operators for type II} and \ref{thm:ladder operators
for type I}, one readily derives differential equations of order
$r + 1$ for the type II and the type I multiple orthogonal
polynomials. As an illustration, we show how to derive differential
equations satisfied by $P_{\vec{n}}(x)$ with $r=2$, the idea of
which can, of course, be extended to general situations. With $r=2$ in
\eqref{def:vec P}, i.e.,
$\textbf{P}_{\vec{n}}(x)=(P_{\vec{n}}(x),P_{\vec{n}-\vec{e}_1}(x),P_{\vec{n}-\vec{e}_2}(x))^T$,
it follows from \eqref{eq:matrix form type II} and straightforward
calculations that
\begin{align*}    %\label{diff eq1}
\textbf{P}_{\vec{n}}'(x)&=N(x)\textbf{P}_{\vec{n}}(x),
\\   %\label{diff eq2}
\textbf{P}_{\vec{n}}''(x)&=N_1(x)\textbf{P}_{\vec{n}}(x),
\\   %\label{diff eq3}
\textbf{P}_{\vec{n}}'''(x)&=N_2(x)\textbf{P}_{\vec{n}}(x),
\end{align*}
where
\begin{align*}
N_1(x)&=N'(x)+N^2(x), \\
N_2(x)&=N''(x)+2N'(x)N(x)+N(x)N'(x)+N^3(x),
\end{align*}
and the derivative of a matrix valued function is understood in the
entry-size manner. Hence,
\begin{align}
P_{\vec{n}}'(x)&=(N)_{00}(x)P_{\vec{n}}(x)+(N)_{01}(x)P_{\vec{n}-\vec{e}_1}(x)+
(N)_{02}(x)P_{\vec{n}-\vec{e}_2}(x), \label{eq:first derivative}
\\
P_{\vec{n}}''(x)&=(N_1)_{00}(x)P_{\vec{n}}(x)+(N_1)_{01}(x)P_{\vec{n}-\vec{e}_1}(x)+
(N_1)_{02}(x)P_{\vec{n}-\vec{e}_2}(x), \label{eq:sec derivative}
\\
P_{\vec{n}}'''(x)&=(N_2)_{00}(x)P_{\vec{n}}(x)+(N_2)_{01}(x)P_{\vec{n}-\vec{e}_1}(x)+
(N_2)_{02}(x)P_{\vec{n}-\vec{e}_2}(x), \label{eq:3rd derivative}
\end{align}
where we use the notation $(M)_{ij}$ to denote the $(i,j)$-th entry
of any given matrix $M$. Then, we can represent
$P_{\vec{n}-\vec{e}_1}$ and $P_{\vec{n}-\vec{e}_2}$ in terms of
$P_{\vec{n}}$, $P_{\vec{n}}'$ and $P_{\vec{n}}''$ by solving
\eqref{eq:first derivative} and \eqref{eq:sec derivative}. Finally,
replacing $P_{\vec{n}-\vec{e}_1}$ and $P_{\vec{n}-\vec{e}_2}$ in
\eqref{eq:3rd derivative} by these relations will lead us to the
third order differential equations satisfied by $P_{\vec{n}}$. Since
the exact formulas are cumbersome, we plan not to write them down
here, but will present some concrete examples in Section
\ref{sec:example}.

\subsection{Compatibility conditions}
We finally state the compatibility conditions for the ladder
equations for multiple orthogonal polynomials.
First we write the nearest-neighbor recurrence relations in a vector form.
Let $\textbf{P}_{\vec{n}}(x)$ be defined by \eqref{def:vec P}, then we
have
\begin{equation}\label{eq:difference relation}
\textbf{P}_{\vec{n}+\vec{e}_l}(x)= W(\vec{n}+\vec{e}_l;x)
\textbf{P}_{\vec{n}}(x),
\end{equation}
for $l=1,\ldots,r$, where
\begin{equation}\label{def:W}
W(\vec{n}+\vec{e}_l;x)=\begin{pmatrix}
    x-b_{\vec{n},l} & -a_{\vec{n},1} & -a_{\vec{n},2} & \cdots & -a_{\vec{n},l} & \cdots & -a_{\vec{n},r} \\
      1 & B_{l,1}(\vec{n}) & 0 & \cdots & 0 & \cdots & 0 \\
      1 &  0 & B_{l,2}(\vec{n})  &  & 0 & \cdots & 0 \\
      \vdots & \vdots & & \ddots &  & \vdots & \vdots \\
      1 & 0 & 0 & \cdots & 0 & \cdots & 0  \\
      \vdots & \vdots & \vdots & & & \ddots  \\
      1 & 0 & 0 & \cdots & 0 & 0 & B_{l,r}(\vec{n})
      \end{pmatrix}
\end{equation}
is a matrix polynomial of degree 1 with
\begin{equation*}
B_{l,j}(\vec{n})=b_{\vec{n}-\vec{e}_j,j}-b_{\vec{n}-\vec{e}_j,l},\qquad
j=1,\ldots,r.
\end{equation*}
The formula \eqref{eq:difference relation} follows from the
nearest-neighbor recurrence relations \eqref{eq:1.10}, since by
eliminating $xP_{\vec{n}}(x)$ between $l$-th and $j$-th relation in
\eqref{eq:1.10}, we have
\begin{equation*}
P_{\vec{n}+\vec{e}_l}(x) - P_{\vec{n}+\vec{e}_j}(x) =
(b_{\vec{n},j}-b_{\vec{n},l}) P_{\vec{n}}(x).
\end{equation*}
Rewriting these relations and the $l$-th relation in \eqref{eq:1.10}
gives us \eqref{eq:difference relation} and \eqref{def:W}. We then
have

\begin{theorem}\label{thm:compa con}
The compatibility conditions
for the ladder equations stated in Theorem \ref{thm:ladder operators for
type II} are given by
\begin{equation}\label{eq:comp cond}
N(\vec{n}+\vec{e}_l;x)W(\vec{n}+\vec{e}_l;x)=W'(\vec{n}+\vec{e}_l;x)
+W(\vec{n}+\vec{e}_l;x)N(\vec{n};x),\qquad l=1,\ldots,r,
\end{equation}
where $N(\vec{n};x)$ is defined in \eqref{eq:matrix form type II}.
\end{theorem}
Clearly, by considering \eqref{eq:comp cond} in an entry-size
manner, one has at least $r(r+1)^2$ equalities, which are much more
complicated than the single weight case. Similar compatibility
conditions for ladder equations associated with the type I multiple
orthogonal polynomials can be obtained easily by using
\eqref{eq:comp cond}. We omit the results here.

\subsection{Outline of the paper}
The rest of this paper is mainly devoted to the proofs of our
theorems. We will prove Theorems \ref{thm:ladder operators for type
II} and \ref{thm:ladder operators for type I} in Section
\ref{sec:proof of 2.1&2.2}. Unlike the treatment of the single
weight case \cite{chen1}, Riemann-Hilbert (RH) problems and the
Christoffel-Darboux formula for multiple orthogonal polynomials will
be two fundamental components in the derivations of the ladder
equations; see also \cite{Gru} for similar treatment to matrix
orthogonal polynomials. We then derive the compatibility conditions
in Section \ref{sec:comp condition}, where we also give a review of
partial difference equations satisfied by the nearest-neighbor
recurrence coefficients, obtained in \cite{WVA}, for later use. We
conclude this paper by applying our results to several concrete
examples, namely, multiple Hermite and Laguerre polynomials, and
multiple orthogonal polynomials with exponential weights and cubic
potentials.

\section{Proofs of Theorems
\ref{thm:ladder operators for type II} and \ref{thm:ladder operators
for type I}} \label{sec:proof of 2.1&2.2}

As mentioned before, the proofs of Theorems \ref{thm:ladder operators
for type II} and \ref{thm:ladder operators for type I} rely on the
Riemann-Hilbert (RH) problem and the Christoffel-Darboux formula
for multiple orthogonal polynomials. In what follows, we first give
a brief introduction of these two aspects.

\subsection{Riemann-Hilbert problems and the Christoffel-Darboux formula for multiple orthogonal polynomials}
The usual orthogonal polynomials can be characterized by a RH
problem of size $2\times 2$ \cite{Fokas}. It was shown by Van Assche
et al. \cite{WVAGerKui} that multiple orthogonal polynomials can
also be described in terms of a RH problem, but for matrices of
order $(r+1)\times(r+1)$.

For convenience, we assume that all the measures are supported on an
oriented, unbounded and simple curve $\Gamma$ in the complex plane.
If the weight functions $w_j$ are H\"{o}lder continuous, we look for
an $(r+1) \times (r+1)$ matrix valued function $Y$ satisfying the
following RH problem:
\begin{enumerate}
\item $Y$ is analytic on $\mathbb{C} \setminus \Gamma$.
\item  $Y$ possesses continuous boundary values
  $Y_{+}(x)$ (from the positive side of $\Gamma$) and $Y_{-}(x)$ (from the negative side of $\Gamma$)
, which satisfies
\begin{equation*} %\label{eq:2.1}
      Y_+(x) = Y_-(x) \begin{pmatrix} 1 &  w_1(x) & w_2(x) & \cdots & w_r(x) \\
                                      0 & 1 & 0 & \cdots & 0 \\
                                      0 & 0 & 1 & & 0 \\
                                     \vdots & \vdots&  & \ddots & 0\\
                                      0 & 0 & \cdots & 0 & 1
                        \end{pmatrix}, \qquad x \in \Gamma.
\end{equation*}
\item As $z\to \infty$, $z\in\mathbb{C}\setminus \Gamma$, we have
\begin{equation*} %\label{eq:2.2}
     Y(z) = \left( I + \O(1/z) \right)\diag(z^{|\vec{n}|},z^{-n_1},z^{-n_2}
\ldots,z^{-n_r}).
\end{equation*}
\end{enumerate}
If the indices $\vec{n}$, $\vec{n}-\vec{e}_j$, $j=1,\ldots,r$, are
normal, then  there exists a unique solution in terms of the type II
multiple orthogonal polynomials which is given by
\begin{equation}  \label{eq:RHPtypeII}
Y=  \begin{pmatrix}
    P_{\vec{n}}(z) & \displaystyle \frac{1}{2\pi i} \int_\Gamma \frac{P_{\vec{n}}(t)w_1(t)}{t-z}\, dt & \cdots &
      \displaystyle    \frac{1}{2\pi i} \int_\Gamma \frac{P_{\vec{n}}(t)w_r(t)}{t-z}\, dt  \\[15pt]
    -2\pi i\gamma_{\vec{n},1} P_{\vec{n}-\vec{e}_1}(z) & \displaystyle
-\gamma_{\vec{n},1} \int_\Gamma
\frac{P_{\vec{n}-\vec{e}_1}(t)w_1(t)}{t-z}\, dt & \cdots &
  \displaystyle -\gamma_{\vec{n},1} \int_\Gamma \frac{P_{\vec{n}-\vec{e}_1}(t)w_r(t)}{t-z}\, dt \\
  \vdots & \vdots & \cdots & \vdots \\
-2\pi i\gamma_{\vec{n},r} P_{\vec{n}-\vec{e}_r}(z) & \displaystyle
-\gamma_{\vec{n},r} \int_\Gamma
\frac{P_{\vec{n}-\vec{e}_r}(t)w_1(t)}{t-z}\, dt & \cdots &
 \displaystyle -\gamma_{\vec{n},r} \int_\Gamma \frac{P_{\vec{n}-\vec{e}_r}(t)w_r(t)}{t-z}\, dt
   \end{pmatrix},
\end{equation}\label{def:gamma k}
where
\begin{equation}\label{2.4} \frac{1}{\gamma_{\vec{n},j}} = \int_\Gamma x^{n_j-1}
P_{\vec{n}-\vec{e}_j}(t)w_j(t)\, dt, \qquad j=1,\ldots,r.
\end{equation}

There exists a similar RH problem for the type I multiple orthogonal
polynomials to determine a  matrix valued function  $X$ of dimension
$(r+1)\times(r+1)$ such that
\begin{enumerate}
\item $X$ is analytic on $\mathbb{C} \setminus \Gamma$.
\item For $x\in\Gamma$, we have
\begin{equation*} %\label{eq:2.1}
      X_+(x) = X_-(x) \begin{pmatrix} 1 &  0 & 0 & \cdots & 0 \\
                                      -w_1(x) & 1 & 0 & \cdots & 0 \\
                                      -w_2(x) & 0 & 1 & & 0 \\
                                     \vdots & \vdots&  & \ddots & 0\\
                                      -w_r(x) & 0 & \cdots & 0 & 1
                        \end{pmatrix}.
\end{equation*}
\item As $z\to \infty$, $z\in\mathbb{C}\setminus \Gamma$, we have
\begin{equation*} %\label{eq:2.2}
     X(z) = \left( I + \O(1/z) \right)\diag(z^{-|\vec{n}|},z^{n_1},z^{n_2}
\ldots,z^{n_r}).
\end{equation*}
\end{enumerate}
Assume that the indices $\vec{n}$, $\vec{n}-\vec{e}_j$,
$j=1,\ldots,r$, are normal, then the unique solution   $X$ is
given by
\begin{equation}\label{eq:X}
  X= \begin{pmatrix}
  \displaystyle \int_\Gamma \frac{Q_{\vec{n}}(t)}{z-t}\, dt
     & 2\pi i A_{\vec{n},1}(z) & \cdots & 2\pi i A_{\vec{n},r}(z) \\
     \displaystyle \frac{c_{\vec{n},1}}{2\pi i} \int_\Gamma \frac{Q_{\vec{n}+\vec{e}_1}(t)}{z-t}\, dt &
     c_{\vec{n},1} A_{\vec{n}+\vec{e}_1,1}(z) & \cdots &
         c_{\vec{n},1} A_{\vec{n}+\vec{e}_1,r}(z) \\
       \vdots & \vdots & \cdots & \vdots \\
     \displaystyle \frac{c_{\vec{n},r}}{2\pi i} \int_\Gamma \frac{Q_{\vec{n}+\vec{e}_r}(t)}{z-t}\, dt &
     c_{\vec{n},r} A_{\vec{n}+\vec{e}_r,1}(z) & \cdots &
         c_{\vec{n},r} A_{\vec{n}+\vec{e}_r,r}(z)
   \end{pmatrix},
\end{equation}
where $Q_{\vec{n}}$ is defined in \eqref{def:Qn} and
\begin{equation}\label{def:cj}
c_{\vec{n},j}=\frac{1}{\kappa_{\vec{n}+\vec{e}_j,j}}, \qquad
j=1,\ldots,r,
\end{equation}
with $\kappa_{\vec{n}+\vec{e}_j,j}$ being the leading coefficient of
$A_{\vec{n}+\vec{e}_j,j}$.%ismail22

If the curve $\Gamma$ is bounded, extra conditions at the endpoints
of $\Gamma$ are required in the RH problem to ensure the uniqueness
of the solution, if the solution exists.

A simple relation between the matrix functions $X$ for the type I
and $Y$ for the type II multiple orthogonal polynomials is given by
Mahler's relation (cf. \cite[Theorem 23.8.3]{Ismail}):
\begin{equation}\label{eq:mahler relation}
X(z) = Y^{-T}(z),
\end{equation}
where $A^{-T}$ is the transpose of the inverse of the matrix $A$.
As we shall see later on, this relation plays an important role in
our proof.

We conclude this section with the Christoffel-Darboux formula for
multiple orthogonal polynomials. This formula was given by Daems and
Kuijlaars \cite{DaeKuij} as stated in the following theorem (see
also \cite{WVA} for another proof):

\begin{theorem}
Suppose $(\vec{n}_i)_{i=0,1,\ldots,|\vec{n}|}$ is a path in
$\mathbb{N}^r$ such that $\vec{n}_0 = \vec{0}$, $\vec{n}_{|\vec{n}|}
= \vec{n}$ and for every $i \in \{0,1,\ldots,|\vec{n}|-1\}$ one has
$\vec{n}_{i+1} - \vec{n}_i = \vec{e}_k$ for some $k \in
\{1,2,\ldots,r\}$. Then
\begin{equation}  \label{CD}
 (x-y) \sum_{i=0}^{|\vec{n}|-1} P_{\vec{n}_i}(x) Q_{\vec{n}_{i+1}}(y)
    = P_{\vec{n}}(x) Q_{\vec{n}}(y) - \sum_{j=1}^r a_{\vec{n},j} P_{\vec{n}-\vec{e}_j}(x)Q_{\vec{n}+\vec{e}_j}(y).
\end{equation}
\end{theorem}

Observe that the right hand side is independent of the path
$(\vec{n}_i)_{i=0,1,\ldots,|\vec{n}|}$ in $\mathbb{N}^r$. A similar
formula for multiple orthogonal polynomials of mixed type is derived
in \cite{AFM}. This formula is of importance in the analysis of
certain random matrices \cite{BleKuij2,Kuijl} and non-intersecting
Brownian motions \cite{DaeKuij2}.

It is worthwhile to point out that if the system of measures is such that
all the indices are normal (perfect system) and the functions
\[  \{ w_1,xw_1,\ldots,x^{n_1-1}w_1,w_2,xw_2,\ldots,x^{n_2-1}w_2,\ldots,
     w_r,xw_r,\ldots,x^{n_r-1}w_r \} \]
are linearly independent (this is exactly our assumption on the $r$
measures in Theorem \ref{thm:ladder operators for type II}), then the
Christoffel-Darboux formula also holds in a component-wise manner,
i.e.,
\begin{equation}\label{eq:CD for component}
(x-y) \sum_{i=0}^{|\vec{n}|-1} P_{\vec{n}_i}(x)
A_{\vec{n}_{i+1},k}(y) = P_{\vec{n}}(x) A_{\vec{n},k}(y) -
\sum_{j=1}^r a_{\vec{n},j}
P_{\vec{n}-\vec{e}_j}(x)A_{\vec{n}+\vec{e}_j,k}(y),
\end{equation}
for $y\in\supp(\mu_k)$ and $k=1,\ldots, r$. To see this, we rewrite
\eqref{CD} as
\begin{equation*}
\sum_{k=1}^{r}R_k(x,y)w_k(y)=0,
\end{equation*}
where
\begin{equation*}
R_{k}(x,y)=(x-y) \sum_{i=0}^{|\vec{n}|-1} P_{\vec{n}_i}(x)
A_{\vec{n}_{i+1},k}(y) - P_{\vec{n}}(x) A_{\vec{n},k}(y) +
\sum_{j=1}^r a_{\vec{n},j}
P_{\vec{n}-\vec{e}_j}(x)A_{\vec{n}+\vec{e}_j,k}(y).
\end{equation*}
With $x$ fixed, it is easily seen that each $R_i(x,y)$ is a
polynomial in $y$ of degree less or equal to $n_k$. The linear
independence then implies $R_i(x,y)=0$, which is \eqref{eq:CD for
component}.

\subsection{Proof of Theorem \ref{thm:ladder operators for type II}}
Let us define an $(r+1)\times(r+1)$ matrix $M$ by
\begin{equation}\label{def:M}
M(x)=M(\vec{n};x)=(M_{ij}(x))_{0\leq i,j\leq r}:=Y'(x)Y^{-1}(x),
\end{equation}
or, equivalently,
\begin{equation}  \label{eq:dy=yM}
Y'(x)=M(x)Y(x),
\end{equation}
where $Y$ is given by \eqref{eq:RHPtypeII}. Since the first column
of $Y$ is given in terms of the type II multiple orthogonal
polynomials, it is easily seen that
\begin{equation}\label{eq:dpn 1}
P_{\vec{n}}'(x)=M_{00}(x)P_{\vec{n}}(x)-\sum_{j=1}^r 2\pi i
\gamma_{\vec{n},j}M_{0j}(x)P_{\vec{n}-\vec{e}_j}(x),
\end{equation}
and
\begin{equation}\label{eq:dpnei 1}
P_{\vec{n}-\vec{e}_i}'(x)=-\frac{1}{2 \pi i
\gamma_{\vec{n},i}}\left(M_{i0}(x)P_{\vec{n}}(x)-\sum_{j=1}^r 2\pi i
\gamma_{\vec{n},j}M_{ij}(x)P_{\vec{n}-\vec{e}_j}(x)\right),
\end{equation}
for $i=1,\ldots,r$.

In the following we use the notation
\begin{equation*}
C(f)(x)=\int_\Gamma\frac{f(t)}{t-x}\,dt,\qquad
x\in\mathbb{C}\setminus\Gamma,
\end{equation*}
for the Cauchy transform of $f$. The derivative $d\, C(f)(x)/dx$
is denoted by $C(f)'(x)$.

On the other hand, the Mahler relation \eqref{eq:mahler relation}
gives
\begin{equation*}
M(x)=Y'(x)X^{T}(x).
\end{equation*}
It then follows from \eqref{eq:RHPtypeII} and \eqref{eq:X} that
\begin{align}\label{eq:M00}
M_{00}(x)&=-P_{\vec{n}}'(x)C(Q_{\vec{n}})(x)+\sum_{k=1}^{r}A_{\vec{n},k}(x)C(P_{\vec{n}}w_k)'(x),\\
M_{0j}(x)&=\frac{c_{\vec{n},j}}{2\pi
i}\left(-P_{\vec{n}}'(x)C(Q_{\vec{n}+\vec{e}_j})(x)
+\sum_{k=1}^{r}A_{\vec{n}+\vec{e}_j,k}(x)C(P_{\vec{n}}w_k)'(x)\right)
\end{align}
with $j=1,\ldots,r$. For $i,j=1,\ldots,r$, we have
\begin{align}\label{eq:Mi0}
M_{i0}(x)&=2 \pi i
\gamma_{\vec{n},i}\left(P_{\vec{n}-\vec{e}_i}'(x)C(Q_{\vec{n}})(x)
-\sum_{k=1}^{r}A_{\vec{n},k}(x)C(P_{\vec{n}-\vec{e}_i}w_k)'(x)\right),\\
M_{ij}(x)&=\gamma_{\vec{n},i}c_{\vec{n},j}\left(P_{\vec{n}-\vec{e}_i}'(x)C(Q_{\vec{n}+\vec{e}_j})(x)
-\sum_{k=1}^{r}A_{\vec{n}+\vec{e}_j,k}(x)C(P_{\vec{n}-\vec{e}_i}w_k)'(x)\right).\label{eq:Mij}
\end{align}
%with $j=1,\ldots,r$.

 The entries $M_{ij}$ in $M$ can be simplified
in view of the following two elementary facts.

\begin{proposition}
 If
\begin{equation*}
\int_\Gamma t^k f(t) \, dt=0,\qquad k=0,\ldots,n-1,
\end{equation*}
then we have
\begin{equation*}
p(x)C(f)(x)=C(pf)(x)
\end{equation*}
for any  polynomial $p$ of degree less than or equal to $n$.
\end{proposition}
Indeed, $$C(pf)(x)=\int_\Gamma
\frac{f(t)(p(t)-p(x))}{t-x}dt+p(x)C(f)(x)$$ and $ (p(t)-p(x))/(t-x)$
is a polynomial of degree less than or equal to $n-1$, hence, by
assumption, the integral disappears.

Integration by parts gives us the following statement.
\begin{proposition} If $f$ is a differentiable and integrable
function that vanishes at the endpoints of $\Gamma$, one has
\begin{equation*}
C(f)'(x)=C(f')(x).
\end{equation*}
\end{proposition}

Hence, the orthogonality of $P$ and $Q$ (see \eqref{eq:orthogonality
of Pn} and \eqref{eq:biorthogonality}) implies
\begin{align}
P_{\vec{n}}'(x)C(Q_{\vec{n}})(x)&=C(P_{\vec{n}}'Q_{\vec{n}})(x), \\
A_{\vec{n},k}(x)C(P_{\vec{n}}w_k)'(x)&
=C(A_{\vec{n},k}P_{\vec{n}}w_k)'(x)-C(P_{\vec{n}}w_k)(x)A_{\vec{n},k}'(x),
\quad k=1,\ldots,r. \label{AC(Pw)}
\end{align}
Inserting the above formulas into \eqref{eq:M00}, it is readily seen
that
\begin{align}\label{eq:M00 simp}
M_{00}(x)&=-C(P_{\vec{n}}'Q_{\vec{n}})(x)
+\sum_{k=1}^{r}\left(C(A_{\vec{n},k}P_{\vec{n}}w_k)'(x)-C(P_{\vec{n}}w_k)(x)
A_{\vec{n},k}'(x)\right)
\nonumber\\
&=-C(P_{\vec{n}}'Q_{\vec{n}})(x)+C(P_{\vec{n}}Q_{\vec{n}})'(x)-\sum_{k=1}^{r}
C(P_{\vec{n}}w_k)(x)A_{\vec{n},k}'(x)
\nonumber \\
&=C(P_{\vec{n}}Q_{\vec{n}}')(x)-\sum_{k=1}^{r}C(P_{\vec{n}}w_kA_{\vec{n},k}')(x)
\nonumber \\
&=C\left(P_{\vec{n}}\sum_{k=1}^r
A_{\vec{n},k}w_k'\right)(x)=-C\left(P_{\vec{n}}\sum_{k=1}^r
A_{\vec{n},k}v_k'w_k\right)(x),
\end{align}
where $v_k(x):=-\ln w_k(x)$. In a similar way, we have
\begin{align}\label{eq:Moj simp}
M_{0j}(x)&=-\frac{c_{\vec{n},j}}{2\pi
i}C\left(P_{\vec{n}}\sum_{k=1}^rA_{\vec{n}+\vec{e}_j,k}v_k'w_k\right)(x),
\\
\label{eq:Mi0 simp} M_{i0}(x)&=2 \pi i \gamma_{\vec{n},i}
C\left(P_{\vec{n}-\vec{e}_i}\sum_{k=1}^{r}A_{\vec{n},k}v_k'w_k\right)(x).
\end{align}

The simplification of $M_{ij}$ with $i,j\geq1$, however, is a
little bit tricky. The problem is that
$A_{\vec{n}+\vec{e}_j,k}(x)C(P_{\vec{n}-\vec{e}_i}w_k)'(x)$ in
\eqref{eq:Mij} can not be written in a form similar to
\eqref{AC(Pw)} (indeed, for $i=j=k$ by comparing the degrees of
$A_{\vec{n}+\vec{e}_j,k}$ and orthogonality relations
(\ref{eq:orthogonality of Pn}), we see that assumptions of the
first proposition are not fulfilled). To overcome this difficulty,
recall that $YY^{-1}=YX^{T}=I$, the $(i,j)$ entry on both sides
gives
\begin{equation*}
\gamma_{\vec{n},i}
c_{\vec{n},j}\left(P_{\vec{n}-\vec{e}_i}(x)C(Q_{\vec{n}+\vec{e}_j})(x)-
\sum_{k=1}^{r}A_{\vec{n}+\vec{e}_j,k}(x)C(P_{\vec{n}-\vec{e}_i}w_k)(x)\right)=\delta_{i,j},\quad
i,j\geq 1.
\end{equation*}
Note that
$P_{\vec{n}-\vec{e}_i}C(Q_{\vec{n}+\vec{e}_j})=C(P_{\vec{n}-\vec{e}_i}Q_{\vec{n}+\vec{e}_j})$.
Taking the derivative of the above formula with respect to $x$, we see
\begin{equation}\label{eq:ij entry}
\begin{aligned}
-\sum_{k=1}^{r}A_{\vec{n}+\vec{e}_j,k}(x)C(P_{\vec{n}-\vec{e}_i}w_k)'(x)
=&\sum_{k=1}^{r}A_{\vec{n}+\vec{e}_j,k}'(x)C(P_{\vec{n}-\vec{e}_i}w_k)(x)
\\
&-C(P_{\vec{n}-\vec{e}_i}'Q_{\vec{n}+\vec{e}_j})(x)
-C(P_{\vec{n}-\vec{e}_i}Q_{\vec{n}+\vec{e}_j}')(x).
\end{aligned}
\end{equation}
%A combination of \eqref{eq:Mij} and \eqref{eq:ij entry} gives us
Using \eqref{eq:ij entry} in \eqref{eq:Mij} we have
\begin{align}
M_{ij}(x)&=\gamma_{\vec{n},i}c_{\vec{n},j}\sum_{k=1}^r\left(
A_{\vec{n}+\vec{e}_j,k}'(x)C(P_{\vec{n}-\vec{e}_i}w_k)(x)
-C(P_{\vec{n}-\vec{e}_i}Q_{\vec{n}+\vec{e}_j}')(x) \right)\nonumber
\\
\label{eq:Mij simp}
&=\gamma_{\vec{n},i}c_{\vec{n},j}
C\left(\sum_{k=1}^rP_{\vec{n}-\vec{e}_i}A_{\vec{n}+\vec{e}_j,k}v_k'w_k\right)(x),
\quad i,j\geq 1.
\end{align}

Replacing $M_{ij}$ in \eqref{eq:dpn 1} and \eqref{eq:dpnei 1} by
\eqref{eq:M00 simp}--\eqref{eq:Mi0 simp} and \eqref{eq:Mij simp}, it
follows
\begin{equation}\label{eq:dpn 2}
\begin{aligned}
P_{\vec{n}}'(x)&=-C\left(P_{\vec{n}}\sum_{k=1}^r
A_{\vec{n},k}v_k'w_k\right)(x)P_{\vec{n}}(x) \\
&~~~+\sum_{j=1}^r
a_{\vec{n},j}C\left(P_{\vec{n}}\sum_{k=1}^rA_{\vec{n}+\vec{e}_j,k}v_k'w_k\right)(x)P_{\vec{n}-\vec{e}_j}(x),
\end{aligned}
\end{equation}
and
\begin{equation}\label{eq:dpnei 2}
\begin{aligned}
P_{\vec{n}-\vec{e}_i}'(x)&=-C\left(P_{\vec{n}-\vec{e}_i}\sum_{k=1}^{r}A_{\vec{n},k}(x)v_k'w_k\right)(x)
P_{\vec{n}}(x)
\\&~~~+\sum_{j=1}^r a_{\vec{n},j}
C\left(\sum_{k=1}^rP_{\vec{n}-\vec{e}_i}A_{\vec{n}+\vec{e}_j,k}v_k'w_k\right)(x)P_{\vec{n}-\vec{e}_j}(x),
\end{aligned}
\end{equation}
for $i=1,\ldots,r$, where we have used   the fact that
\begin{equation}\label{eq:a=cgamma}
c_{\vec{n},j}\gamma_{\vec{n},j}=a_{\vec{n},j};
\end{equation}
see \eqref{eq:an representation}, \eqref{2.4}, \eqref{eq:kappa} and %\eqref{def:gamma k} and
\eqref{def:cj}. %{\textbf{Here maybe (2.4) instead of (2.1) and
%maybe also we should add (1.12).}}

To show \eqref{lowering II}, we observe from \eqref{eq:dpn 2} that
\begin{align}\label{eq:dpn 3}
P_{\vec{n}}'(x)=&~P_{\vec{n}}(x)\int_\Gamma P_{\vec{n}}(t)
\sum_{k=1}^r A_{\vec{n},k}(t)\frac{v_k'(t)-v_k'(x)}{x-t} w_k(t)\,dt
\nonumber
\\ &-\sum_{j=1}^r a_{\vec{n},j} P_{\vec{n}-\vec{e}_j}(x)\int_\Gamma
P_{\vec{n}}(t)  \sum_{k=1}^r
A_{\vec{n}+\vec{e}_j,k}(t)\frac{v_k'(t)-v_k'(x)}{x-t}
w_k(t)\,dt+E(x),
\end{align}
where
\begin{equation}\label{def:E}
E(x)=\sum_{k=1}^r \left(\sum_{j=1}^r a_{\vec{n},j}C(
P_{\vec{n}}A_{\vec{n}+\vec{e}_j,k}w_k)(x)P_{\vec{n}-\vec{e}_j}(x)
-C( P_{\vec{n}}A_{\vec{n},k}w_k)(x)P_{\vec{n}}(x)\right)v_{k}'(x).
\end{equation}
It is then equivalent to show $E=0$. To this end, we apply the
Christoffel-Darboux formula \eqref{eq:CD for component} to each
coefficient of $v_k'$ in \eqref{def:E} and obtain
\begin{align*} %\label{eq:coefficeint of vk}
&\sum_{j=1}^r a_{\vec{n},j}C(
P_{\vec{n}}A_{\vec{n}+\vec{e}_j,k}w_k)(x)P_{\vec{n}-\vec{e}_j}(x)
-C( P_{\vec{n}}A_{\vec{n},k}w_k)(x)P_{\vec{n}}(x) \nonumber \\
&=\sum_{i=0}^{|\vec{n}|-1}\left(\int_\Gamma P_{\vec{n}}(t)
A_{\vec{n}_{i+1},k}(t)w_k(t) \,dt\right) P_{\vec{n}_i}(x)=0,
\end{align*}
where we have used  the orthogonality of $P_{\vec{n}}$
\eqref{eq:orthogonality of Pn} in the last equality. Hence, it is
immediate that $E=0$ and \eqref{lowering II} follows.

In a similar way, we see from \eqref{eq:dpnei 2} that
\begin{align}\label{eq:dpnei 3}
P_{\vec{n}-\vec{e}_i}'(x)=&~P_{\vec{n}}(x)\int_\Gamma
P_{\vec{n}-\vec{e}_i}(t)\sum_{k=1}^r
A_{\vec{n},k}(t)\frac{v_k'(t)-v_k'(x)}{x-t}w_k(t)\,dt\nonumber\\
&-\sum_{j=1}^r a_{\vec{n},j} P_{\vec{n}-\vec{e}_j}(x)\int_\Gamma
P_{\vec{n}-\vec{e}_i}(t)  \sum_{k=1}^r
A_{\vec{n}+\vec{e}_j,k}(t)\frac{v_k'(t)-v_k'(x)}{x-t}w_k(t)\,dt
\nonumber
\\&+E_i(x),
\end{align}
where
\begin{equation*} %\label{def:Ei}
E_i(x)=\sum_{k=1}^r \left(\sum_{j=1}^r a_{\vec{n},j}C(
P_{\vec{n}-\vec{e}_i}A_{\vec{n}+\vec{e}_j,k}w_k)(x)P_{\vec{n}-\vec{e}_j}(x)
-C(P_{\vec{n}-\vec{e}_i}A_{\vec{n},k}w_k)(x)P_{\vec{n}}(x)\right)v_{k}'(x).
\end{equation*}
Again with the aid of the Christoffel-Darboux formula \eqref{eq:CD for
component} and the orthogonality of $P_{\vec{n}-\vec{e}_i}$, it is
readily  checked that
\begin{align*}
&\sum_{j=1}^r a_{\vec{n},j}C(
P_{\vec{n}-\vec{e}_i}A_{\vec{n}+\vec{e}_j,k}w_k)(x)P_{\vec{n}-\vec{e}_j}(x)
-C(P_{\vec{n}-\vec{e}_i}A_{\vec{n},k}w_k)(x)P_{\vec{n}}(x) \nonumber \\
&=\sum_{i=0}^{|\vec{n}|-1}\left(\int_\Gamma P_{\vec{n}-\vec{e}_i}(t)
A_{\vec{n}_{i+1},k}(t)w_k(t) \,dt\right)
P_{\vec{n}_i}(x)=P_{\vec{n}-\vec{e}_i}(x)\delta_{i,j}.
\end{align*}
This, together with \eqref{eq:dpn 3}, implies
\begin{equation}\label{eq:Ei explicit}
E_i(x)=P_{\vec{n}-\vec{e}_i}(x)v_i'(x).
\end{equation}
Combining \eqref{eq:dpnei 3} and \eqref{eq:Ei explicit} then leads
to \eqref{raising II}.
This completes the proof of Theorem \ref{thm:ladder operators for
type II}.
\medskip

\noindent \textbf{Remark.}
There is a simple relation between the matrices $M$ in \eqref{def:M}
and $N$ in \eqref{def:N}. Note that the first column on both sides
of \eqref{eq:dy=yM} gives us
\begin{equation*}  %\label{N eq}
\frac{d}{dx}
\begin{pmatrix} P_{\vec{n}}(x) \\
-2 \pi i \gamma_{\vec{n},1}P_{\vec{n}-\vec{e}_1}(x)\\
\vdots \\
-2 \pi i \gamma_{\vec{n},r}P_{\vec{n}-\vec{e}_r}(x)
\end{pmatrix}
=M(x)\begin{pmatrix} P_{\vec{n}}(x) \\
-2 \pi i \gamma_{\vec{n},1}P_{\vec{n}-\vec{e}_1}(x)\\
\vdots \\
-2 \pi i \gamma_{\vec{n},r}P_{\vec{n}-\vec{e}_r}(x)
\end{pmatrix}.
\end{equation*}
We then obtain from \eqref{eq:matrix form type II} that
\begin{equation}\label{eq:M and N}
N(x)=\diag\left(1,-\frac{1}{2 \pi i
\gamma_{\vec{n},1}},\cdots,-\frac{1}{2 \pi i
\gamma_{\vec{n},r}}\right)M(x) \diag\left(1,-2 \pi i
\gamma_{\vec{n},1},\cdots,-2 \pi i \gamma_{\vec{n},r}\right).
\end{equation}

\subsection{Proof of Theorem \ref{thm:ladder operators for type I}}
We shall prove Theorem \ref{thm:ladder operators for type I} by
establishing \eqref{eq:matrix form type I}. The formulas
\eqref{raising I} and \eqref{lowering I} then follow from
straightforward calculations using \eqref{eq:matrix form type I}
and \eqref{def:N}.

To see \eqref{eq:matrix form type I}, we first take a derivative
with respect to $x$ on both sides of the equation $Y(x)Y^{-1}(x)=I$
and obtain from \eqref{def:M} that
\begin{equation*}
Y(x)(Y^{-1})'(x)=-M(x).
\end{equation*}
Equivalently, in view of \eqref{eq:mahler relation},
\begin{equation}\label{eq:DXT =-XTM}
(X^T)'(x)=-X^T(x)M(x).
\end{equation}
Next, we observe from \eqref{eq:X} that, by transposing the $l$-th
($l=1,\cdots,r$) row on both sides of \eqref{eq:DXT =-XTM}, we get
\begin{equation*}
\frac{d}{dx}
\begin{pmatrix} 2 \pi i A_{\vec{n},l}(x) \\
c_{\vec{n},1}A_{\vec{n}+\vec{e}_1,l}(x)\\
\vdots \\
c_{\vec{n},r}A_{\vec{n}+\vec{e}_r,l}(x)
\end{pmatrix}
=-M^T(x)\begin{pmatrix} 2 \pi i A_{\vec{n},l}(x) \\
c_{\vec{n},1}A_{\vec{n}+\vec{e}_1,l}(x)\\
\vdots \\
c_{\vec{n},r}A_{\vec{n}+\vec{e}_r,l}(x)
\end{pmatrix},
\end{equation*}
that is,
\begin{align*}
\frac{d}{dx}
\begin{pmatrix}  A_{\vec{n},l}(x) \\
A_{\vec{n}+\vec{e}_1,l}(x)\\
\vdots \\
A_{\vec{n}+\vec{e}_r,l}(x)
\end{pmatrix}
=&-\diag\left(\frac{1}{2 \pi i},
\frac{1}{c_{\vec{n},1}},\cdots,\frac{1}{c_{\vec{n},r}}
\right)M^T(x)\diag\left(2 \pi i, c_{\vec{n},1}, \cdots,
c_{\vec{n},r}\right) \nonumber \\
&\times
\begin{pmatrix} A_{\vec{n},l}(x) \\
A_{\vec{n}+\vec{e}_1,l}(x)\\
\vdots \\
A_{\vec{n}+\vec{e}_r,l}(x)
\end{pmatrix}.
\end{align*}
Replacing $M^T$ in the above formula with the aid of \eqref{eq:M and
N}, we finally arrive at \eqref{eq:matrix form type I} after simple
manipulations using \eqref{eq:a=cgamma}.
This completes the proof of Theorem \ref{thm:ladder operators for
type I}.

%---------------------------------------------------------------------
\section{Proof of Theorem \ref{thm:compa con}}
\label{sec:comp condition}

To show \eqref{eq:comp cond}, we note that equations
\eqref{eq:matrix form type II} and \eqref{eq:difference relation}
can be viewed as a Lax pair. Differentiating both sides of
\eqref{eq:difference relation}, together with \eqref{eq:matrix form
type II}, implies
\begin{equation*}
N(\vec{n}+\vec{e}_l;x)\textbf{P}_{\vec{n}+\vec{e}_l}(x)=\left(W'(\vec{n}+\vec{e}_l;x)
+W(\vec{n}+\vec{e}_l;x)N(\vec{n};x)\right)\textbf{P}_{\vec{n}}(x),
\end{equation*}
for $l=1,\ldots,r$. Using \eqref{eq:difference relation} again, it
follows
\begin{equation} \label{eq:comp cond with P}
N(\vec{n}+\vec{e}_l;x)W(\vec{n}+\vec{e}_l;x)\textbf{P}_{\vec{n}}(x)=\left(W'(\vec{n}+\vec{e}_l;x)
+W(\vec{n}+\vec{e}_l;x)N(\vec{n};x)\right)\textbf{P}_{\vec{n}}(x).
\end{equation}
It is now sufficient to show that the components of the vector
$\textbf{P}_{\vec{n}}$ are linearly independent whenever $\vec{n}$
is a normal index. This assertion is already shown in \cite{WVA}.
For the convenience of the reader, we repeat the argument here.
Suppose that $(c_0,c_1,\ldots,c_r)$
are such that
\[  c_0 P_{\vec{n}} + \sum_{j=1}^r c_j P_{\vec{n}-\vec{e}_j} = 0, \]
then by comparing the leading coefficients, it follows $c_0=0$.
If we multiply by $x^{n_k-1}$ and integrate with respect to $\mu_k$, then
\[      c_k \int x^{n_k-1} P_{\vec{n}-\vec{e}_k}(x)\, d\mu_k(x) = 0.  \]
We claim that the above integral does not vanish. Indeed, if the
integral vanishes, then $P_{\vec{n}} - a P_{\vec{n}-\vec{e}_k}$ is a
monic polynomial of degree $|\vec{n}|$ satisfying the orthogonality
relations \eqref{eq:orthogonality of Pn} for any $a$, which contradicts
the normality of $\vec{n}$. As a consequence, $c_k=0$ for $k=1,\ldots,r$.
The linear independence of the polynomials in $\textbf{P}_{\vec{n}}$ and
\eqref{eq:comp cond with P} gives the compatibility conditions \eqref{eq:comp cond}.
This completes the proof of Theorem \ref{thm:compa con}.
\medskip

\noindent \textbf{Remark.}
From \eqref{eq:difference relation}, it is readily seen that one can obtain the vector
$\textbf{P}_{\vec{n}+\vec{e}_i+\vec{e}_j}$ in two different ways. The first way is
to first compute $\textbf{P}_{\vec{n}+\vec{e}_i}$ from $\textbf{P}_{\vec{n}}$
and then to compute $\textbf{P}_{\vec{n}+\vec{e}_i+\vec{e}_j}$ from
$\textbf{P}_{\vec{n}+\vec{e}_i}$:
\begin{equation*}
\textbf{P}_{\vec{n}+\vec{e}_i+\vec{e}_j}(x)= W(\vec{n}+\vec{e}_i+\vec{e}_j;x)
W(\vec{n}+\vec{e}_i;x)\textbf{P}_{\vec{n}}(x).
\end{equation*}
The second way is
to first compute $\textbf{P}_{\vec{n}+\vec{e}_j}$ from $\textbf{P}_{\vec{n}}$
and then to compute $\textbf{P}_{\vec{n}+\vec{e}_i+\vec{e}_j}$
from $\textbf{P}_{\vec{n}+\vec{e}_j}$:
\begin{equation*}
\textbf{P}_{\vec{n}+\vec{e}_i+\vec{e}_j}(x)
=W(\vec{n}+\vec{e}_j+\vec{e}_i;x)W(\vec{n}+\vec{e}_j;x)\textbf{P}_{\vec{n}}(x).
\end{equation*}
Hence, we obtain another kind of compatibility condition
\begin{equation*}
W(\vec{n}+\vec{e}_i+\vec{e}_j;x)W(\vec{n}+\vec{e}_i;x)=W(\vec{n}+\vec{e}_j+\vec{e}_i;x)
W(\vec{n}+\vec{e}_j;x).
\end{equation*}
This was observed in \cite{WVA}, and from the explicit formula of
$W$ the following nearest-neighbor recurrence relations for multiple
orthogonal polynomials were obtained.
\begin{theorem}\label{thm:gen compatibility}
\textup{\cite[Theorem 3.2]{WVA}} Suppose all multi-indices $\vec{n} \in
\mathbb{N}^r$ are normal. Suppose $1 \leq i \neq j \leq r$, then the
recurrence coefficients for the nearest neighbor recurrence
relations (\ref{eq:1.10}) satisfy
\begin{align*}
b_{\vec{n}+\vec{e}_i,j} - b_{\vec{n},j} &= b_{\vec{n}+\vec{e}_j,i} - b_{\vec{n},i},     %\label{eq:r1}
\\
\sum_{k=1}^r a_{\vec{n}+\vec{e}_j,k} - \sum_{k=1}^r a_{\vec{n}+\vec{e}_i,k}
       &= \det \begin{pmatrix}  b_{\vec{n}+\vec{e}_j,i} & b_{\vec{n},i} \\
b_{\vec{n}+\vec{e}_i,j} & b_{\vec{n},j} \end{pmatrix}, %\label{eq:r2}
\\
\frac{a_{\vec{n},i}}{a_{\vec{n}+\vec{e}_j,i}} &=
\frac{b_{\vec{n}-\vec{e}_i,j}-b_{\vec{n}-\vec{e}_i,i}}{b_{\vec{n},j}-b_{\vec{n},i}}.
%\label{eq:r3}
\end{align*}
\end{theorem}
This theorem tells us that the recurrence coefficients of multiple orthogonal
polynomials cannot be arbitrary, but satisfy a system of partial
difference equations. This is not the case for the usual orthogonal
polynomials.

\section{Some examples}\label{sec:example}

In this section, we shall derive the ladder equations and differential equations for
several examples of multiple orthogonal polynomials. For convenience, it is assumed
that $r=2$ throughout this section. The (normal) multi-index $\vec{n}$ is now given by
$(n,m)\in\mathbb{N}^2$. We also use the following notation for the recurrence coefficients:
\begin{equation*}   %\label{correspondence}
a_{\vec{n},1}=a_{n,m},\;\;a_{\vec{n},2}=b_{n,m},\;\;b_{\vec{n},1}=c_{n,m},\;\;b_{\vec{n},2}=d_{n,m}.
\end{equation*}
The recurrence relations \eqref{eq:1.10} then read
\begin{align}%\label{rec rel r=2}
    xP_{n,m}(x) &= P_{n+1,m}(x)+c_{n,m}P_{n,m}(x)
    +a_{n,m} P_{n-1,m}(x) + b_{n,m} P_{n,m-1}(x), \label{eq:2.6} \\
    xP_{n,m}(x) &= P_{n,m+1}(x)+d_{n,m}P_{n,m}(x)
    + a_{n,m} P_{n-1,m}(x) + b_{n,m} P_{n,m-1}(x), \label{eq:2.7}
\end{align}
with $a_{0,m}=0$ and $b_{n,0}=0$ for all $n,m \geq 0$. In view of
Theorem \ref{thm:gen compatibility}, the following relations hold:
\begin{align}
    d_{n+1,m}-d_{n,m} &= c_{n,m+1} - c_{n,m},  \label{eq:3.1}  \\
    b_{n+1,m}-b_{n,m+1} + a_{n+1,m}-a_{n,m+1} &=
    \det \begin{pmatrix} d_{n+1,m} & d_{n,m} \\
    c_{n,m+1} & c_{n,m} \end{pmatrix}, \label{eq:3.2} \\
    \frac{a_{n,m+1}}{a_{n,m}} &= \frac{c_{n,m}-d_{n,m}}{c_{n-1,m}-d_{n-1,m}},  \label{eq:3.3} \\
    \frac{b_{n+1,m}}{b_{n,m}} &= \frac{c_{n,m}-d_{n,m}}{c_{n,m-1}-d_{n,m-1}}.  \label{eq:3.4}
\end{align}

\subsection{Multiple Hermite polynomials}
Multiple Hermite polynomials $H_{n,m}$ are type II multiple
orthogonal polynomials defined by
\begin{align*}
\int_{-\infty}^\infty x^k H_{n,m}(x) e^{-x^2+c_1x}\,dx &=0,
\qquad k = 0, 1, \ldots, n-1, \\
\int_{-\infty}^\infty x^k H_{n,m}(x) e^{-x^2+c_2x}\,dx &=0, \qquad k
= 0, 1, \ldots, m-1,
\end{align*}
where $c_1\neq c_2$; cf. \cite{BleKuij}, \cite[\S 23.5]{Ismail} and
\cite[\S 3.4]{WVAEC}. Hence, the functions $v_k(x)$ in Theorem
\ref{thm:ladder operators for type II} are equal to $x^2-c_k x$. A
simple calculation with the aid of \eqref{eq:orthogonality of Pn}
then gives
\begin{equation}   \label{eq:N in MHP}
N(n,m;x)=\begin{pmatrix}
0 & 2 a_{n,m} & 2 b_{n,m} \\
-2 & 2x-c_1 & 0 \\
-2 & 0 & 2x-c_2
\end{pmatrix}.
\end{equation}
Recall that $N(n,m;x)$ is defined by \eqref{eq:matrix form type II}
and explicitly given in \eqref{def:Nij}. This, together with
\eqref{def:W} and the compatibility conditions \eqref{eq:comp cond}, leads to the
the following system of difference equations for the recurrence coefficients:
\begin{align}
2a_{n,m}-2a_{n+1,m}+2b_{n,m}-2b_{n+1,m}+1&=0,\label{eq:5.8}\\
2a_{n,m}-2a_{n,m+1}+2b_{n,m}-2b_{n,m+1}+1&=0,\label{eq:5.9}\\
2c_{n,m}=c_1,\;\;2d_{n,m}&=c_2,\label{eq:c(n,m) & d(n,m)}\\
2(d_{n,m-1}-c_{n,m-1})-(c_2-2c_{n,m})&=0,\label{eq:5.11}\\
2b_{n+1,m}(d_{n,m-1}-c_{n,m-1})-b_{n,m}(c_2-2c_{n,m})&=0,\label{eq:5.12} \\
2(d_{n-1,m}-c_{n-1,m})+c_1-2d_{n,m}&=0,\label{eq:5.13}\\
2a_{n,m+1}(d_{n-1,m}-c_{n-1,m})+a_{n,m}(c_1-2d_{n,m})&=0\label{eq:5.14}.
\end{align}
It turns out one can easily obtain the recurrence coefficients explicitly
from these equations. Indeed, from \eqref{eq:c(n,m) & d(n,m)},
it follows
\begin{equation*}
c_{n,m}=c_1/2,\qquad d_{n,m}=c_2/2.
\end{equation*}
Since $c_1\neq c_2$, equations \eqref{eq:5.11} and \eqref{eq:5.12} imply that
$b_{n,m}$ is independent of $n$, while \eqref{eq:5.13} and \eqref{eq:5.14} show
that $a_{n,m}$ is independent of $m$.
Hence, by \eqref{eq:5.8} and \eqref{eq:5.9}, it follows
\begin{equation}\label{eq:mid step for a and b}
2a_{n,m}-2a_{n,m+1}+1=0,\quad 2b_{n,m}-2b_{n+1,m}+1=0.
\end{equation}
Note that our initial conditions are $a_{0,m}=0$ and $b_{n,0}=0$, we then
obtain from \eqref{eq:mid step for a and b}
\begin{equation*}
a_{n,m}=n/2, \qquad b_{n,m}=m/2.
\end{equation*}
Thus, we arrive at the following lowering and raising equations for
multiple Hermite polynomials $H_{n,m}$:
\begin{align}\label{eq:a}
H_{n,m}'(x)&=n H_{n-1,m}(x)+m H_{n,m-1}(x),
\\
\label{eq:b}
H_{n-1,m}'(x)&=-2H_{n,m}(x)+(2x-c_1)H_{n-1,m}(x),
\\
\label{eq:c}
H_{n,m-1}'(x)&=-2H_{n,m}(x)+(2x-c_2)H_{n,m-1}(x).
\end{align}
Oberve that equation (\ref{eq:b}) (similarly (\ref{eq:c}))
can be alternatively written as
$$\left(e^{-x^2+c_1 x}H_{n-1,m}(x)\right)'=-2e^{-x^2+c_1 x}H_{n,m}(x).$$
These formulas are not new but can already be found in \cite[\S 23.8.2]{Ismail}.

%Some of these equations can be easily solved. Moreover, the recurrence coefficients for multiple Hermite polynomials are known explicitly \cite{Ismail, WVA}, in particular,
%\begin{gather*}
%a_{n,m}=n/2,\;\;b_{n,m}=m/2,\;\;c_{n,m}=c_1/2,\;\;d_{n,m}=c_2/2.
%\end{gather*}
%Moreover, calculating explicitly,
%$$c_1(n,m)=-\frac{i 2^{m+n}(c1-c2)^{-m} e^{-c1^2/4}\sqrt{\pi}}{(n-1)!},\;\;c_2(n,m)=-\frac{i 2^{m+n}(c2-c1)^{-n} e^{-c2^2/4}\sqrt{\pi}}{(m-1)!}.$$
%One can check that these expressions satisfy the  system of difference equations derived above.

Finally, by \eqref{eq:N in MHP} and arguments at the end of Section
\ref{sec:ladder equations}, we see that the type II multiple
orthogonal polynomial $H_{n,m}$ satisfies the following linear
differential equation of order 3:
\begin{multline*}
p'''(x)+ (c_1+c_2-4x)p''(x)+\left(c_1(c_2-2x)+2(m+n-1-c_2 x+2x^2)\right)p'(x)\\
+2(c_1 m+c_2 n-2(m+n)x)p(x)=0.
\end{multline*}
Similarly, we can derive a third order differential equation for the
type I multiple Hermite polynomials $A_{(n,m),l},\;l=1,2$, which is
given by
\begin{multline*}
q'''(x)-(c_1+c_2-4x)q''(x)+\left(c_1(c_2-2x)+2(m+n-1-c_2 x+2x^2)\right)q'(x)\\
-2(c_1 m+c_2 n-2(m+n)x)q(x)=0.
\end{multline*}
We point out that the above differential equation is independent of
$l$, which can also be seen from \eqref{eq:matrix form type I}.

%----------------------------------------------------------------------------------
\subsection{Multiple Laguerre polynomials of the second kind}\
\label{sec:MLP 2nd kind}
These polynomials are defined by the
orthogonality conditions
\begin{align*}
\int_0^\infty x^k L_{n,m}(x) x^{\alpha} e^{-c_1 x}\,dx &=0,
\qquad k = 0, 1, \ldots, n-1, \\
\int_0^\infty x^k L_{n,m}(x) x^{\alpha} e^{-c_2 x}\,dx &=0, \qquad k
= 0, 1, \ldots, m-1,
\end{align*}
where we assume that $\alpha>0$ and $c_1,c_2>0$ with $c_1\neq c_2$;
cf. \cite{BleKuij}, \cite[Remark 5 on p.~160]{NikSor}, \cite[\S
23.4.2]{Ismail} and \cite[\S 3.3]{WVAEC}. Hence, the functions
$v_k(x)$ in Theorem \ref{thm:ladder operators for type II} are equal
to $-\alpha \ln x+c_k x$, $k=1,2$. An appeal to the biorthogonality
conditions \eqref{eq:biorthogonality} implies that
\begin{equation}   \label{eq:N in MLP}
N(n,m;x)= \frac{1}{x}\begin{pmatrix}
-\alpha\int_0^{\infty}\frac{L_{n,m}(t)Q_{n,m}(t)}{t}\,dt & A(n,m)
& B(n,m) \\
-\alpha\int_0^{\infty}\frac{L_{n-1,m}(t)Q_{n,m}(t)}{t}\,dt & C(n,m)-\alpha+c_1x& D(n,m) \\
-\alpha\int_0^{\infty}\frac{L_{n,m-1}(t)Q_{n,m}(t)}{t}\,dt & E(n,m)
& F(n,m)-\alpha+c_2x
\end{pmatrix},
\end{equation}
where $Q_{n,m}$ is defined in \eqref{def:Qn} with
$w_j(x)=x^{\alpha}e^{-c_j x}$, $j=1,2$, and
\begin{align}
A(n,m)&=\alpha
       a_{n,m}\int_0^{\infty}\frac{L_{n,m}(t)Q_{n+1,m}(t)}{t}\,dt,\label{eq:A(n,m)}\\
B(n,m)&=\alpha
       b_{n,m}\int_0^{\infty}\frac{L_{n,m}(t)Q_{n,m+1}(t)}{t}\,dt,\label{eq:B(n,m)}\\
C(n,m)&=\alpha
       a_{n,m}\int_0^{\infty}\frac{L_{n-1,m}(t)Q_{n+1,m}(t)}{t}\,dt,\label{eq:C(n,m)}\\
D(n,m)&=\alpha
       b_{n,m}\int_0^{\infty}\frac{L_{n-1,m}(t)Q_{n,m+1}(t)}{t}\,dt,\label{eq:D(n,m)}\\
E(n,m)&=\alpha
       a_{n,m}\int_0^{\infty}\frac{L_{n,m-1}(t)Q_{n+1,m}(t)}{t}\,dt,\label{eq:E(n,m)}\\
F(n,m)&=\alpha
       b_{n,m}\int_0^{\infty}\frac{L_{n,m-1}(t)Q_{n,m+1}(t)}{t}\,dt,\label{eq:F(n,m)}
\end{align}
are certain constants depending on $n$ and $m$. We want to give an
explicit representation of $N(n,m;x)$ in \eqref{eq:N in MLP}. We
first observe that the first column in \eqref{eq:N in MLP} can be
evaluated by comparing the leading coefficients on both sides of the
differential equations \eqref{lowering II} and \eqref{raising II}.
For instance, by \eqref{lowering II} and \eqref{eq:N in MLP}, it
follows
\begin{multline*}
xL_{n,m}'(x)=
\left(-\alpha\int_0^{\infty}\frac{L_{n,m}(t)Q_{n,m}(t)}{t}dt\right)L_{n,m}(x)\\
+A(n,m)L_{n-1,m}(x)+B(n,m)L_{n,m-1}(x).
\end{multline*}
Recall that $L_{n,m}(x)=x^{n+m}+\ldots$, hence by comparing the
coefficients of order $n+m$ on both sides of the above equation, it
is readily seen that
\begin{equation}\label{eq:first entry}
\alpha\int_0^{\infty}\frac{L_{n,m}(t)Q_{n,m}(t)}{t}dt=-(n+m).
\end{equation}
Similarly, we obtain
\begin{align}
\alpha\int_0^{\infty}\frac{L_{n-1,m}(t)Q_{n,m}(t)}{t}\,dt=c_1,\quad
\alpha\int_0^{\infty}\frac{L_{n,m-1}(t)Q_{n,m}(t)}{t}\,dt=c_2.
\end{align}
To estimate the integrals in \eqref{eq:A(n,m)}--\eqref{eq:F(n,m)},
we shall make use of the compatibility conditions \eqref{eq:comp
cond} and the known results of the recurrence coefficients:
\begin{equation}
a_{n,m}=\frac{(n+m+\alpha)n}{c_1^2},\qquad
b_{n,m}=\frac{(n+m+\alpha)m}{c_2^2},
\end{equation}
and
\begin{equation}\label{eq:c d in MLP}
c_{n,m}=\frac{2n+m+\alpha+1}{c_1}+\frac{m}{c_2},\qquad
d_{n,m}=\frac{n+2m+\alpha+1}{c_2}+\frac{n}{c_1};
\end{equation}
see \cite[\S 5.4]{WVA}. The compatibility conditions \eqref{eq:comp
cond} in this case are given by
\begin{multline}\label{eq:comp cond 1 MLP}
x N(n+1,m;x)
\begin{pmatrix} x-c_{n,m} & -a_{n,m}
& -b_{n,m} \\
1 & 0 & 0 \\
1 & 0 & d_{n,m-1}-c_{n,m-1}
\end{pmatrix}
\\=\begin{pmatrix} x & 0
& 0 \\
0 & 0 & 0 \\
0 & 0 & 0
\end{pmatrix}
+\begin{pmatrix} x-c_{n,m} & -a_{n,m}
& -b_{n,m} \\
1 & 0 & 0 \\
1 & 0 & d_{n,m-1}-c_{n,m-1}
\end{pmatrix}x N(n,m;x),
\end{multline}
and
\begin{multline}\label{eq:comp cond 2 MLP}
x N(n,m+1;x)
\begin{pmatrix} x-d_{n,m} & -a_{n,m}
& -b_{n,m} \\
1 & c_{n-1,m}-d_{n-1,m} & 0 \\
1 & 0 & 0
\end{pmatrix}
\\=\begin{pmatrix} x & 0
& 0 \\
0 & 0 & 0 \\
0 & 0 & 0
\end{pmatrix}
+\begin{pmatrix} x-d_{n,m} & -a_{n,m}
& -b_{n,m} \\
1 & c_{n-1,m}-d_{n-1,m} & 0 \\
1 & 0 & 0
\end{pmatrix}x N(n,m;x).
\end{multline}
With the aid of \eqref{eq:N in MLP} and the explicit formulas
\eqref{eq:first entry}--\eqref{eq:c d in MLP}, we see from the
$(2,2)$-entry in \eqref{eq:comp cond 1 MLP} that
\begin{equation*}
A(n,m)=c_1a_{n,m}=\frac{(n+m+\alpha)n}{c_1}.
\end{equation*}
Similarly, the $(3,3)$-entry in \eqref{eq:comp cond 2 MLP} implies
\begin{equation*}
B(n,m)=c_2b_{n,m}=\frac{(n+m+\alpha)m}{c_2}.
\end{equation*}
Then, we observe from the $(2,3)$-entry in \eqref{eq:comp cond 1 MLP}
that
\begin{equation*}
c_1b_{n,m}+(d_{n,m-1}-c_{n,m-1})D(n+1,m)=B(n,m),
\end{equation*}
which gives
\begin{equation*}
D(n,m)=-\frac{mc_1}{c_2}.
\end{equation*}
Similarly, the $(3,2)$-entry in \eqref{eq:comp cond 2 MLP} implies that
\begin{equation*}
E(n,m)=-\frac{nc_2}{c_1},
\end{equation*}
which is independent of $m$. Finally, we obtain from the
$(2,1)$-entry in \eqref{eq:comp cond 1 MLP} that
\begin{equation*}
c_1c_{n,m}+C(n+1,m)-\alpha+D(n+1,m)=n+m,
\end{equation*}
hence,
\begin{equation*}
C(n,m)=-n.
\end{equation*}
An appeal to the $(3,1)$-entry in \eqref{eq:comp cond 1 MLP} gives
\begin{equation*}
c_2c_{n,m}+E(n+1,m)+F(n+1,m)-\alpha=n+m+c_2(c_{n,m-1}-d_{n,m-1}),
\end{equation*}
thus,
\begin{equation*}
F(n,m)=-m.
\end{equation*}
Combining all these results, we find the following lowering and
raising equations for multiple Laguerre polynomials of the second
kind $L_{n,m}$:
\begin{align*}
xL_{n,m}'(x)&=(n+m)L_{n,m}(x)+\left(\frac{(n+m+\alpha)n}{c_1}\right)
L_{n-1,m}(x) \nonumber \\
&~~~~~~~~~~~~~~~~~~~
~~~~~+\left(\frac{(n+m+\alpha)m}{c_2}\right)L_{n,m-1}(x),
\\
xL_{n-1,m}'(x)&=-c_1L_{n,m}(x)-\left(n+\alpha-c_1x\right)
L_{n-1,m}(x)-\frac{mc_1}{c_2}L_{n,m-1}(x),
\\
xL_{n,m-1}'(x)&=-c_2L_{n,m}(x)-\frac{nc_2}{c_1}L_{n-1,m}(x)
-\left(m+\alpha-c_2x\right)L_{n,m-1}(x).
\end{align*}
Furthermore, the third order differential equation satisfied by
$L_{n,m}$ is given by
\begin{align*}
x^2p'''(x)&-\left(x^2(c_1+c_2)-2x(\alpha+1)\right)p''(x) \nonumber
\\&+\left(x^2c_1c_2-x[(c_1+c_2)(\alpha+1)
-nc_1-mc_2]+\alpha(\alpha+1)\right)p'(x) \nonumber
\\
&-\left(xc_1c_2(n+m)-\alpha(nc_1+mc_2)\right)p(x)=0,
\end{align*}
which can also be found in \cite[\S 4.3]{AptBraWVA}. Similarly, it
can be shown that the type I multiple Laguerre polynomials of the
second kind satisfy the following differential equation:
\begin{align*}
x^2q'''(x)&+\left(x^2(c_1+c_2)-2x(\alpha-1)\right)q''(x) \nonumber
\\&+\left(x^2c_1c_2-x[(c_1+c_2)(\alpha-1)
-nc_1-mc_2]+\alpha(\alpha-1)\right)q'(x) \nonumber
\\
&+\left(xc_1c_2(n+m)-\alpha(nc_1+mc_2)\right)q(x)=0.
\end{align*}
%----------------------------------------------------------------------------------
\subsection{Multiple Laguerre polynomials of the first kind}\
These polynomials are defined by the orthogonality conditions
\begin{align*}
\int_0^\infty x^k L_{n,m}(x) x^{\alpha_1} e^{-x}\,dx &=0,
\qquad k = 0, 1, \ldots, n-1, \\
\int_0^\infty x^k L_{n,m}(x) x^{\alpha_2} e^{-x}\,dx &=0, \qquad k =
0, 1, \ldots, m-1,
\end{align*}
where $\alpha_1,\alpha_2>0$ and $\alpha_1-\alpha_2\notin
\mathbb{Z}$; cf. \cite{BleKuij}, \cite[\S 23.4.1]{Ismail} and
\cite[\S 3.2]{WVAEC}. The functions $v_k(x)$ in Theorem
\ref{thm:ladder operators for type II} are now equal to $-\alpha_k
\ln x+ x$, $k=1,2$. A straightforward calculation using
\eqref{eq:biorthogonality} gives
\begin{equation*}   %\label{eq:N in MLP 1st kind}
N(n,m;x)= \frac{1}{x}\begin{pmatrix}
-\int_0^{\infty}\frac{L_{n,m}(t)Q_{n,m}^*(t)}{t}\,dt & A^*(n,m)
& B^*(n,m) \\
-\int_0^{\infty}\frac{L_{n-1,m}(t)Q_{n,m}^*(t)}{t}\,dt & C^*(n,m)-\alpha_1+x& D^*(n,m) \\
-\int_0^{\infty}\frac{L_{n,m-1}(t)Q_{n,m}^*(t)}{t}\,dt & E^*(n,m) &
F^*(n,m)-\alpha_2+x
\end{pmatrix},
\end{equation*}
where we define
\begin{equation*}
Q_{n,m}^*(t)=\left(\alpha_1A_{(n,m),1}(t)t^{\alpha_1}+
\alpha_2A_{(n,m),2}(t)t^{\alpha_2}\right)e^{-t},
\end{equation*}
with $A_{(n,m),j}$, $j=1,2$ the associated type I multiple
orthogonal polynomials and
\begin{align*}
A^*(n,m)&=  a_{n,m}\int_0^{\infty}\frac{L_{n,m}(t)Q_{n+1,m}^*(t)}{t}\,dt,  %\label{eq:A(n,m) MLP1}
\\
B^*(n,m)&=  b_{n,m}\int_0^{\infty}\frac{L_{n,m}(t)Q_{n,m+1}^*(t)}{t}\,dt, %\label{eq:B(n,m)MLP1}
\\
C^*(n,m)&=  a_{n,m}\int_0^{\infty}\frac{L_{n-1,m}(t)Q_{n+1,m}^*(t)}{t}\,dt, %\label{eq:C(n,m)MLP1}
\\
D^*(n,m)&=  b_{n,m}\int_0^{\infty}\frac{L_{n-1,m}(t)Q_{n,m+1}^*(t)}{t}\,dt,  %\label{eq:D(n,m)MLP1}
\\
E^*(n,m)&=  a_{n,m}\int_0^{\infty}\frac{L_{n,m-1}(t)Q_{n+1,m}^*(t)}{t}\,dt,  %\label{eq:E(n,m)MLP1}
\\
F^*(n,m)&=  b_{n,m}\int_0^{\infty}\frac{L_{n,m-1}(t)Q_{n,m+1}^*(t)}{t}\,dt,  %\label{eq:F(n,m)MLP1}
\end{align*}
are certain constants depending on $n$ and $m$. Again by comparing
the leading coefficients on both sides of the differential equations
\eqref{lowering II} and \eqref{raising II}, we have
\begin{align*}
\int_0^{\infty}\frac{L_{n,m}(t)Q_{n,m}^*(t)}{t}\,dt&=-(n+m) \\
\int_0^{\infty}\frac{L_{n-1,m}(t)Q_{n,m}^*(t)}{t}\,dt&=
\int_0^{\infty}\frac{L_{n,m-1}(t)Q_{n,m}^*(t)}{t}\,dt=1.
\end{align*}
Since the recurrence coefficients are explicitly given by
\begin{align*}
a_{n,m}&=\frac{n(n+\alpha_1)(n+\alpha_1-\alpha_2)}{n-m+\alpha_1-\alpha_2},
&&\qquad
b_{n,m}=\frac{m(m+\alpha_2)(m+\alpha_2-\alpha_1)}{m-n+\alpha_2-\alpha_1},\\
c_{n,m}&=2n+m+\alpha_1+1,&&\qquad d_{n,m}=n+2m+\alpha_2+1,
\end{align*}
see \cite[\S 5.3]{WVA}, the same strategy as in Section \ref{sec:MLP
2nd kind} gives
\begin{align*}
A^*(n,m)&=\frac{n(n+\alpha_1)(n+\alpha_1-\alpha_2)}{n-m+\alpha_1-\alpha_2},
&&\quad
B^*(n,m)=\frac{m(m+\alpha_2)(m+\alpha_2-\alpha_1)}{m-n+\alpha_2-\alpha_1},\\
C^*(n,m)&=-n, &&\quad F^*(n,m)=-m,
\\D^*(n,m)&=E^*(n,m)=0.
\end{align*}
Hence, we find the following ladder equations for multiple Laguerre
polynomials of the first kind $L_{n,m}$:
\begin{align*}
xL_{n,m}'(x)&=(n+m)L_{n,m}(x)+
\left(\frac{n(n+\alpha_1)(n+\alpha_1-\alpha_2)}{n-m+\alpha_1-\alpha_2}\right)
L_{n-1,m}(x) \nonumber\\
&~~~~~~~~~~~~~~~~~~~
~~~~~~+\left(\frac{m(m+\alpha_2)(m+\alpha_2-\alpha_1)}{m-n+\alpha_2-\alpha_1}\right)L_{n,m-1}(x),
\\
xL_{n-1,m}'(x)&=-L_{n,m}(x)-\left(n+\alpha_1-x\right) L_{n-1,m}(x),
\\
xL_{n,m-1}'(x)&=-L_{n,m}(x)-\left(m+\alpha_2-x\right)L_{n,m-1}(x).
\end{align*}
Furthermore, the third order differential equation satisfied by
$L_{n,m}$ is given by
\begin{align*}
x^2p'''(x)&+\left(-2x^2+(\alpha_1+\alpha_2+3)x\right)p''(x)
\nonumber
\\&+\left(x^2-x(\alpha_1+\alpha_2-n-m+3)+(\alpha_1+1)(\alpha_2+1)\right)p'(x) \nonumber
\\
&-\left(x(n+m)-(n+m+nm+\alpha_1
m+\alpha_2n)\right)p(x)=0,
\end{align*}
which can also be found in \cite[\S 4.3]{AptBraWVA}. Similarly, we
have that the type I multiple Laguerre polynomials of the first kind
satisfy the following differential equation:
\begin{align*}
x^2q'''(x)&+\left(2x^2-(\alpha_1+\alpha_2-3)x\right)q''(x)
\nonumber
\\&+\left(x^2-x(\alpha_1+\alpha_2-n-m-3)+(\alpha_1-1)(\alpha_2-1)\right)q'(x) \nonumber
\\
&+\left(x(n+m)-(mn-n-m+\alpha_1 m+\alpha_2n)\right)q(x)=0.
\end{align*}

%-------------------------------------------------------------------------------------
\subsection{Multiple exponential polynomials with cubic potentials}
These polynomials are defined by the orthogonality conditions
\begin{align*}
\int_\Gamma x^k P_{n,m}(x) e^{-x^3/3-c_1 x}\,dx &=0,
\qquad k = 0, 1, \ldots, n-1, \\
\int_\Gamma x^k P_{n,m}(x) e^{-x^3/3-c_2 x}\,dx &=0, \qquad k = 0,
1, \ldots, m-1,
\end{align*}
where $c_1\neq c_2$ and the contour $\Gamma$ is taken to be the set
$\{x: \Re~x^3<0\}$, or simply, $\{x: \arg x=\pm 2\pi/3\}$. This
definition generalizes the orthogonal polynomials introduced by
Magnus in \cite[\S 6]{magnus}. Since $v_k(x)=x^3/3+c_k x$, we obtain
from (\ref{eq:biorthogonality}) and (\ref{eq:1.10}) (or
(\ref{eq:2.6}) and (\ref{eq:2.7}) in case $r=2$) that
\begin{equation*}
N(n,m;x)=\begin{pmatrix}
-a_{n,m}-b_{n,m} & a_{n,m}(x+c_{n,m}) & b_{n,m}(x+d_{n,m})\\
-x-c_{n-1,m} & x^2+a_{n,m}+c_1&b_{n,m}\\
-x-d_{n,m-1}&a_{n,m}&x^2+b_{n,m}+c_2
\end{pmatrix}.
\end{equation*}
%The matrix $M(x)$ in (\ref{compatibility1}), (\ref{compatibility2}) is given explicitly by
%\begin{equation*}
%M(x)=\begin{pmatrix}
%-a_{n,m}-b_{n,m} & a_{n,m}(x+c_{n,m})/c_1(n,m) & b_{n,m}(x+d_{n,m})/c_2(n,m)\\
%(-x-c_{n-1,m})c_1(n,m) & x^2+a_{n,m}+c_1&b_{n,m}c_1(n,m)/c_2(n,m)\\
%(-x-d_{n,m})c_2(n,m)&a_{n,m}c_2(n,m)/c_1(n,m)&x^2+b_{n,m}+c_2
%\end{pmatrix}.
%\end{equation*}
If we introduce the matrix form for $N$:
$$N(n,m;x)=F_2(n,m)x^2+F_1(n,m)x+F_0(n,m)$$
and similarly for $W$ in \eqref{def:W}:
\begin{equation*}
W(n+1,m;x)=Rx+R_{1},\qquad
W(n,m+1;x)=Rx+R_{2},
\end{equation*}
where
\begin{equation*}
R=\begin{pmatrix} 1 & 0
& 0 \\
0 & 0 & 0 \\
0 & 0 & 0
\end{pmatrix},\quad
\end{equation*}
and
\begin{equation*}
R_{1}=\begin{pmatrix} -c_{n,m} & -a_{n,m}
& -b_{n,m} \\
1 & 0 & 0 \\
1 & 0 & d_{n,m-1}-c_{n,m-1}
\end{pmatrix}, \quad R_2=\begin{pmatrix} -d_{n,m} & -a_{n,m}
& -b_{n,m} \\
1 & c_{n-1,m}-d_{n-1,m} & 0 \\
1 & 0 & 0
\end{pmatrix},
\end{equation*}
then the compatibility conditions \eqref{eq:comp cond} are given by
the following equations (which are not identically zero), after
collecting the coefficients with respect to $x$:
\begin{align*}
F_1(n+1,m)R_{1}+F_0(n+1,m)R&=RF_0(n,m)+R_{1}F_1(n,m),\\
F_1(n,m+1)R_{2}+F_0(n,m+1)R&=RF_0(n,m)+R_{2}F_1(n,m),\\
F_0(n+1,m)R_{1}&=R+R_{1}F_0(n,m),\\
F_0(n,m+1)R_{2}&=R+R_{2}F_0(n,m).
\end{align*}
Although in this case we do not know the exact expressions of the
recurrence coefficients, the above relations indeed give us many
nonzero equations for these coefficients. It can be shown that all
but the following four equations can be simplified using them and
(\ref{eq:3.1}) and (\ref{eq:3.2}):
\begin{align}
c_1+a_{n,m}+a_{n+1,m}+b_{n,m}+b_{n+1,m}+c_{n,m}^2&=0,\label{eq:nonlinear1}\\
c_2+a_{n,m}+a_{n,m+1}+b_{n,m}+b_{n,m+1}+d_{n,m}^2&=0,\label{eq:nonlinear2}\\
\nonumber
-1-b_{n,m}c_{n,m}+b_{n+1,m}c_{n,m}-a_{n,m}(c_{n-1,m}+c_{n,m})\\
+a_{n+1,m}(c_{n,m}+c_{n+1,m})-b_{n,m}d_{n,m-1}+b_{n+1,m}d_{n+1,m}&=0,\\\nonumber
-1-b_{n,m}d_{n,m-1}-b_{n,m}d_{n,m}-a_{n,m}(c_{n-1,m}+d_{n,m})\\
+a_{n,m+1}(c_{n,m+1}+d_{n,m})+b_{n,m+1}d_{n,m+1}+b_{n,m+1}d_{n,m}&=0.
\end{align}
In addition to equations (\ref{eq:3.1})--(\ref{eq:3.4}) one can get,
for instance,
$$\det \begin{pmatrix} d_{n+1,m} & d_{n,m} \\
c_{n,m+1} & c_{n,m} \end{pmatrix} =d_{n,m}^2-c_{n,m}^2-c_1+c_2.$$ In
particular, we see that equations \eqref{eq:nonlinear1} and
\eqref{eq:nonlinear2} are nonlinear, which is different from the
previous examples. We hope these relations will be helpful in
the further study of multiple exponential polynomials with cubic
potentials.

Finally, we point out that it is also possible to derive the
differential equations for the associated type I and type II
multiple orthogonal polynomials. Since they  are cubersome, we shall
not write them down here, but mention that they are of the form
$$a_2(x)p'''(x)+a_4(x)p''(x)+a_6(x)p'(x)+a_5(x)p(x)=0,$$
where $a_k(x)$ are polynomials in $x$ of degree $k$ with the
coefficients depending on the recurrence coefficients
$a_{n,m},b_{n,m},c_{n,m},d_{n,m}$.

%\section{Discussion}
%
%
%In this paper we have derived the ladder operators for multiple
%orthogonal polynomials of the first and the second type. They can be
%used to derive the differential equations for MOPs and the
%compatibility conditions with the recurrence relations. The
%recurrence coefficients of the  recurrence relations satisfy a
%system of nonlinear difference equations. It might be interesting to
%further study  multiple Laguerre, Jacobi and other MOPs, including
%discrete cases and to derive the associated Lax pairs. Moreover, the
%recurrence coefficients satisfy a system of partial differential
%equations  with respect to the parameters in the weight
%\cite{TodaWVA} which can be combining with    the discrete system to
%give new interesting nonlinear differential  systems.
% $x^{\alpha}e^{-x}$, Lax pair will have
%singularities, discrete

\section*{Acknowledgements}

GF is supported by the MNiSzW Iuventus Plus grant Nr
0124/IP3/2011/71. WVA is supported by KU Leuven research grant
OT/08/033 and FWO research grant G.0427.09. LZ is a Postdoctoral
Fellow of the Fund for Scientific Research - Flanders (FWO),
Belgium.

%\begin{verbatim}
%Galina Filipuk
%Faculty of Mathematics, Informatics and Mechanics
%University of Warsaw
%Banacha 2 02-097 Warsaw
%Poland
%filipuk@mimuw.edu.pl
%\end{verbatim}
%
%\begin{verbatim}
%Walter Van Assche
%Department of Mathematics
%Katholieke Universiteit Leuven
%Celestijnenlaan 200B box 2400
%BE-3001 Leuven
%Belgium
%Walter.VanAssche@wis.kuleuven.be
%\end{verbatim}
%
%\begin{verbatim}
%Lun Zhang
%Department of Mathematics
%Katholieke Universiteit Leuven
%Celestijnenlaan 200 B box 2400
%BE-3001 Leuven
%Belgium
%lun.zhang@wis.kuleuven.be
%\end{verbatim}

\end{document}